\documentclass[11pt,a4paper,reqno]{amsart}
\usepackage{tikz}
\usetikzlibrary{arrows.meta}
\usepackage[utf8]{inputenc}
\usepackage{amsmath,bm,bbm, amssymb}
\usepackage{mathtools}
\usepackage{fullpage}
\usepackage{color}
\usepackage{tcolorbox}
\usepackage{float}
\usepackage{subcaption}
\usepackage{animate}
\usepackage{etoolbox}
\usepackage{comment}
\usepackage{pgf}
\usepackage{hyperref}
\usetikzlibrary{calc}
\usetikzlibrary{patterns}
\usetikzlibrary{arrows}
\usetikzlibrary{quotes,angles}
\usetikzlibrary{decorations.pathreplacing}
\usepackage{tkz-euclide}
\usepackage[utf8]{inputenc}
\usepackage{dsfont}
\usepackage{wrapfig}

\def\ms{\mathsf}
\def\mc{\mathcal}

\def\mbb{\mathbb}
\def\mb{\mathbf}

\def\Z{\mbb Z}
\def\E{\mbb E}
\def\EE{\mc E}

\def\P{\mbb P}

\def\R{\mbb R}
\def\N{\mbb N}
\def\LL{\mathcal L}
\def\RR{\mathcal R}
\def\BB{\mc B}
\def\NN{\mathcal N}
\def\SS{\mathcal S}

\def\f{\frac}
\def\co{\colon}

\def\t{\tau}
\def\a{\alpha}
\def\b{\beta}
\def\g{\gamma}
\def\de{\delta}

\def\e{\varepsilon}
\def\eps{\e}
\def\es{\emptyset}
\def\one{\mathbbmss{1}}

\def\k{\kappa}
\def\wh{\widehat}

\def\ff{\infty}
\def\bel{\begin{lemma}}
\def\enl{\end{lemma}}
\def\bepr{\begin{proposition}}
\def\enpr{\end{proposition}}
\def\bep{\begin{proof}}
\def\enp{\end{proof}}

\def\sm{\setminus}

\newcommand{\pimsta}[1]{\pi_{-s, t, a}}

\def\vp{\varphi}
\def\to{\uparrow}

\def\th{\theta}

\def\su{\subseteq}
\def\De{\Delta}

\def\been{\begin{enumerate}}
\def\enen{\end{enumerate}}

	\def\bec{\begin{corollary}}
	\def\enc{\end{corollary}}

\def\bet{\begin{theorem}}
\def\ent{\end{theorem}}

\DeclareMathOperator*{\argmax}{arg\,max}

\DeclareMathOperator*{\diameter}{diam}
\DeclareMathOperator*{\dist}{dist}

\theoremstyle{plain}
\newtheorem{theorem}{Theorem}
\newtheorem{proposition}[theorem]{Proposition}
\newtheorem{corollary}[theorem]{Corollary}
\newtheorem{lemma}{Lemma}

\theoremstyle{definition}

\theoremstyle{remark}
\newtheorem{remark}{Remark}

\usepackage{enumitem}
\usepackage{etoolbox}
	\patchcmd{\section}{\scshape}{\bfseries}{}{}
	\makeatletter
	\renewcommand{\@secnumfont}{\bfseries}
	\makeatother
	
	\usepackage[normalem]{ulem}	
	\usepackage[section]{placeins}


\begin{document}

\author[Christian Hirsch]{Christian Hirsch}
\author[Daniel Willhalm]{Daniel Willhalm}

\address[Christian Hirsch]{Department of Mathematics, Aarhus University, Ny Munkegade 118, 8000 Aarhus C, Denmark}
\email[Christian Hirsch]{hirsch@math.au.dk}
\address[Daniel Willhalm]{Bernoulli Institute, University of Groningen, Nijenborgh 9, 9747 AG Groningen, Netherlands \newline \hspace*{3.147cm}CogniGron (Groningen Cognitive Systems and Materials Center), Nijenborgh 4, 9747 AG Groningen, Netherlands}
\email[Daniel Willhalm]{d.willhalm@rug.nl}

\title{Upper large deviations for power-weighted edge lengths in spatial random networks}
\date{\today}

\begin{abstract} 
We study the large-volume asymptotics of the sum of power-weighted edge lengths $\sum_{e \in E}|e|^\alpha$ in Poisson-based spatial random networks. In the regime $\alpha > d$, we provide a set of sufficient conditions under which the upper large deviations asymptotics are characterized by a condensation phenomenon, meaning that the excess is caused by a negligible portion of Poisson points. Moreover, the rate function can be expressed through a concrete optimization problem. This framework encompasses in particular directed, bidirected and undirected variants of the $k$-nearest neighbor graph, as well as suitable $\beta$-skeletons.

\vspace{4mm}
			
			
\end{abstract}
		
\maketitle

\section{Introduction} \label{intro_sec} Many real-world networks are not merely a collection of nodes and edges but live in an ambient Euclidean space. Thanks to seminal research efforts on laws of large numbers and central limit theorems,  we now have good understanding of how characteristics computed from stochastic models for geometric networks behave on average in large sampling windows,  and how they fluctuate around the mean \cite{weakLLN, normGeomProb}. However,  when envisioning such models to be used in security-critical applications,  it is essential to understand also the behavior in rare events. The theory of large deviations is designed to deal with such questions. Its achievement is to reduce the understanding of rare events to solving deterministic optimization problems.

On a very general level,  one can think of two radically different causes for a rare event that we refer to as \emph{homogenization} and \emph{condensation}, respectively. In the case of homogenization small but consistent deviations throughout the sampling window add up to yield a macroscopic deviation of the considered quantity. On the other hand,  in the case of condensation,  there is a small isolated structure with the property that its configuration is so extraordinary that it is alone responsible for a deviation that is visible on the macroscopic level. We stress that condensation effects are not by any means restricted to spatial random networks but also play an important role in Erd\H{o}s-R\'enyi graphs,  branching processes, mathematical biology and statistical physics  \cite{adams, andreis, betzDereichCondensation, dereichMoertersCondensation, dereich}. In the classical setting of sums of random variables,  this effect is typical for heavy-tailed models.

For network functionals with finite exponential moments, which includes the power-weighted edge lengths for a wide range of graphs in the case that the power is strictly smaller than the dimension, the homogenization could be made rigorous under very general near-additivity and stabilization conditions  \cite{schreiberyukich, yukLDP2}. However,  on the side of condensation, the research is far less well-developed. Recently,  a breakthrough could be achieved by describing the large deviations of seeing too many edges in the Gilbert graph  \cite{chatterjeeharel} based on a Poisson point process in  $\R^d $. Loosely speaking,  these additional edges are induced by a clique obtained from putting a large number of points in a small spatial domain.

In this work,  we illustrate that condensation phenomena in upper large deviations are not restricted to the Gilbert graph but occur for a broad class of spatial random networks,  including most prominently the  $k $-nearest neighbor graph (kNN). To that end,  we study the upper large deviations of the sum of power-weighted edge lengths,  i.e.,   $\sum_e |e|^\a $,  where the sum is taken over all network edges in a growing sampling window and $\a$ denotes the power considered. This is a fundamental characteristic for spatial random networks,  which has already been studied in detail for the Gilbert graph and the directed spanning forest \cite{chinmoy, reitznerschultethale}.

Speaking of the kNN,  for  $k = 1 $ and very large  $\a $,  the excess weight is induced by a single large edge. Although this is no longer the case for general  $k\ge 1 $ and $\a > d$,  we show that the condensate can still be described in terms of a specific spatial optimization problem. Besides $k $-nearest neighbor graphs,  our framework also encompasses circle-based $\b $-skeletons in two dimensions.

The proof idea is to adapt and refine a three-step strategy that has already been successfully implemented to understand the onset of condensation phenomena in other contexts  \cite{sphere, chatterjeeharel}. First,  the proportion of nodes making a very large contribution to the power-weighted edge lengths is negligible. We identify these nodes as the condensate. Second,  the contributions from nodes outside of the condensate sharply concentrates around the mean. Finally,  analyzing the most likely way that the condensate can cause the excess weight leads to the spatial optimization problem mentioned earlier.

The rest of the article is organized as follows. Section \ref{sec_model} contains precise statements of and conditions for our main results on the upper large deviations of the power-weighted edge lengths. Here,  we also describe the spatial optimization problem in detail that determines the shape of the condensate.  In Sections \ref{sec_application} and \ref{sec_application_condensation}, the theorems connecting the upper large deviations to the optimization problem are applied to the directed,  bidirected and undirected version of the kNN as well as two-dimensional circle-based  $\beta $-skeletons for  $\beta>1 $. Lastly,  Sections \ref{sec_proof_1} and \ref{sec_proof_2} deal with the proofs of our results.

\section{Model and main results}  \label{sec_model} 

To assist the reader, we start by loosely collecting some of the most important notation here. Let $d\ge1$ be the dimension. By $|x|$ we denote the Euclidean norm of $x\in\R^d$. For $e=(x, y)\in(\R^d)^2$, we set $|e| \coloneqq |x-y|$, which is interpreted as length of an edge between $x$ and $y$. Given three points $x,y,z\in\R^d$, we denote the absolute value of the angle of the triangle spanned by $x$, $y$ and $z$ at point $y$ by $\angle xyz$. Further, $B_r (x) \coloneqq \{y \in\R^d \co |y - x|\le r\}  $ denotes the Euclidean ball with radius  $r>0 $, centered at $x \in \R^d $ and for a Borel set $C\subseteq\R^d$ we will use $|C|$ to denote the $d$-dimensional Lebesgue measure of $C$. The symbol $\partial$ refers to the boundary operator that can be applied to a subset of $\R^d$. The ceiling function $\lceil\cdot\rceil$ and floor function $\lfloor\cdot\rfloor$ will appear and are given by $\lceil t\rceil \coloneqq \min\{m\in\Z\colon m\ge t\}$ and $\lfloor t\rfloor \coloneqq \max\{m\in\Z\colon m\le t\}$ for $t\in\R$. By $\mathbf N$ and $\mathbf N_0$, we denote the space of all locally finite subsets of $\R^d$, where the latter must additionally contain the origin $0 \in \R^d $. For a configuration $\varphi\in\mathbf{N}$ and a set $C\subseteq\R^d$, by $\varphi(C)$, we mean $\#(\varphi\cap C)$, the number of points in $\varphi$ that are within $C$. Throughout the paper $Q_n \coloneqq [-n/2,n/2]^d$, $n\ge 1$ represents a cubical observation window.

In the following we describe the general graphs that we study. For $\varphi\in\mathbf N$, the pair $G(\varphi) \coloneqq (\varphi, E)$ represents a directed graph, along with a set of edges $E\coloneqq E(\varphi)\subseteq\{(x,y)\colon x\neq y\in\varphi\}$ on the vertex set $\varphi$. In particular, we stress that the edges are drawn according to some general construction rule that does not depend on the specific point configuration and the edge set is determined once we fix $\varphi$ and does not require any randomness. For $\varphi\in\mathbf N_0$, we let
\begin{equation}\label{equation_out_neighbors}
\EE (\varphi) \coloneqq \{z\in \varphi \co  (0,z)\in E (\varphi)\}
\end{equation}
denote the set of out-neighbors of the origin and \begin{equation}\label{equation_in_out_neighbors}
E_0 (\varphi) \coloneqq\EE (\varphi) \cup \{ x \in \varphi\co 0 \in \EE (\varphi-x) + x\}
\end{equation}
all out- and in-neighbors of $0$. Whenever convenient, we use $\EE_x(\psi) \coloneqq \EE(\psi-x) + x$ for the out-neighbors of $x\in\psi\in\mathbf N$ instead.

In this work, we study the upper large deviations of the sum of  $\a $-power-weighted edge lengths in the box  $Q_n $ for  $\a > d $. For $\varphi\in\mathbf N$, that is the quantity 
$$
H_{n,\ms{dir}} ^{ (\a)}  (G (\varphi)) \coloneqq \f1{n^d}  \sum_{\substack{e =  (x, y)\in E\\ x\in \varphi \cap Q_n} }  |e|^\a = \f1{n^d} \sum_{\substack{z \in\EE (\varphi - x)\\x \in \varphi \cap Q_n }}  |z|^\a.
$$
Hence, by defining the score function $\xi^{ (\a)} _\ms{dir}  (\psi) \coloneqq\sum_{z\in\EE (\psi)}  |z|^\a  $ for $\psi \in\mathbf{N}_0$, we can also express  $H_{n, \ms{dir}} ^{ (\a)} (G (\varphi)) $ as 
 \begin{equation}\label{equation_representation_general_limit}
H_{n,\ms{dir}} ^{ (\a)}  (G (\vp)) = \f1{n^d} \sum_{x \in \vp \cap Q_n}   \xi^{ (\a)} _\ms{dir}  (\vp- x).
\end{equation}
If we represent the nodes of a directed graph by a Poisson point process $X \subseteq \R^d $ with intensity  $1 $, then, $G(X)$ plugged into the representation in \eqref{equation_representation_general_limit} embeds our problem in the setting of general limit results in stochastic geometry, where  a score is assigned to each  $x\in X $ encoding the contribution to the total power-weighted edge lengths.  

Moreover, we note that a directed graph naturally gives rise to two further spatial networks, namely an  \emph{undirected network}, where an edge is put between two nodes  $x, y $ if there is a directed edge from   $x $ to  $y $  \emph{or} a directed edge from  $y $ to  $x $, and a  \emph{bidirected network}, where an edge is put between  $x, y $ if there is a directed edge from   $x $ to  $y $  \emph{and}  a directed edge from  $y $ to  $x $, see \cite[Section 2.3]{weakLLN}.  To extend our results also for these networks, we henceforth work with a score function $\xi^{ (\a)}  $ that, for $\varphi\in\mathbf N_0$, may take one of the following three forms
 $$
\xi^{ (\a)}  (\varphi) \coloneqq \begin{cases}  \xi^{ (\a)} _\ms{dir}  (\varphi)\coloneqq \sum_{x\in\EE (\varphi)}  |x|^\a; \\
\xi^{ (\a)} _\ms{undir}  (\varphi) \coloneqq\sum_{x\in\EE (\varphi)}  \f12 |x|^\a + \f12|x|^\a\one\{0\not\in\EE (\varphi-x)\} ;\\
\xi^{ (\a)} _\ms{bidir}  (\varphi) \coloneqq \sum_{x\in\EE (\varphi)}  \f12|x|^\a\one\{0\in\EE (\varphi-x)\}. 
\end{cases} 
 $$
 In words, the definition of $\xi_\ms{undir}^{(\a)}$ means that if $x$ is an out-neighbor of $0$ but not an in-neighbor, then the edge length $|x|$ contributes fully to the score at $0$, whereas it is not considered for the score at $x$.

We proceed by denoting the corresponding functional for $\varphi\in\mathbf{N}$ by
\begin{equation}\label{equation_representation_functional}
H_n^{ (\a)}  (\varphi) \coloneqq\f1{n^d}  \sum_{x\in \varphi\cap Q_n}  \xi^{ (\a)}  (\varphi-x)
\end{equation}
and if we plug the Poisson point process $X$ as random point configuration in \eqref{equation_representation_functional}, we abbreviate
\begin{equation}\label{equation_representation_functional_X}
H_n  \coloneqq H_n^{ (\a)}  (X).
\end{equation}
In order to describe the large-deviation asymptotic for the upper tails of $H_n $, we require that the graph and the score function satisfy some additional properties.  Our conditions are designed having in mind the (undirected/bidirected)  kNN and a version of the  $\beta $-skeleton as prototypical examples, see Section \ref{sec_application}.  It will become apparent that some of the conditions are substantially more delicate than the ones appearing for weak laws of large number  (WLLNs)  or central limit theorems (CLTs) on Poisson functionals \cite{CLT,weakLLN}.  This is because for many of the spatial random networks  satisfying WLLNs and CLTs like Delaunay tessellations  (DTs), Gabriel graphs  (GGs) or relative neighborhood graphs  (RNGs), the upper large deviations will be markedly different from the ones of the kNN. In all of these graphs, the excess in the large deviation tail might be determined by configurations with a growing number of nodes. For instance, the DT, GG and RNG can, with significantly high probability, exhibit a large total sum of power-weighted edge lengths by having more than a negligible proportion of edges almost parallel to each other. Nevertheless, we decided to present our results in a general framework for two reasons.  First, we can pinpoint precisely to the requirements that are not satisfied by standard examples mentioned earlier.  Second, if one aims to establish upper-large deviation asymptotics for a specific class of networks, the conditions give a clear view at which points additional arguments will be needed to prove the desired result. 

We now first state the conditions rigorously.  After that, we include a detailed discussion to explain more precisely their meaning and impact.  We have not attempted to aggressively minimize the number of conditions because this compactification would entail the risk of leading to statements that are less accessible.   The conditions are the following.

\begin{enumerate}  
\item  $\EE $ is  \emph{scale invariant}: $\tau\EE (\varphi)=\EE (\tau\varphi) $ for all $\varphi\in\mathbf N_0 $ and  $\tau>0 $.
\item Adding a new point  affects only a bounded number of nodes: there exists $c_\ms{FIN} > 0$ such that for every  $y\in\R^d $ and  $\varphi\in\mathbf N $,\begin{equation} \label{FIN} \tag{\textbf{FIN}} \#\big\{x\in \varphi \co \EE (\varphi-x) \neq\EE ( (\varphi-x)\cup\{y-x\} )\big\}  \le c_\ms{FIN}.   \end{equation}  
\item  $\EE $ has \emph{bounded large edge density}: there exists  $c_\ms{FIN2} \geq 1 $ such that for all  $M>0 $ and  $\varphi\in\mathbf N $,\begin{equation} \label{FIN2} \tag{\textbf{FIN2}} \#\big\{x\in \varphi\cap B_M (0)\co\max_{y\in\EE (\varphi-x)}  |y| >M\big\}   \le c_\ms{FIN2}. \end{equation}  
	\item Proceeding in the vein of \cite{weakLLN}, we introduce a stabilization condition for  $G $. This condition is based on a collection of cones $S_i $,  $i\le I_d$ with apex 0 whose union covers the whole space and which do not have parts of their lateral boundary parallel to any coordinate axis of $\R^d$. Then, for a constant $c_\ms{STA} > 0$ and $\vp \in \mb N_0$, we put 
$$\SS_i (\varphi) \coloneqq c_\ms{STA}\inf\{r>0\co\varphi (S_i\cap B_r (0)) \ge c_\ms{STA} \}.$$
		We say that  $G $ is \emph{stabilizing}  if there exists  $c_\ms{STA} \ge 1 $ such that for every $\eta\in\mathbf N_0 $ there exists  $\mathbf N_0\ni\theta\subseteq\eta $ such that   (i)  $\theta \su\cup_{i \le I_d}  \big (S_i\cap B_{\SS_i (\eta)}  (0)\big) \eqqcolon\mathcal{B} $, (ii) $\#\theta\le I_d c_\ms{STA}$, and (iii)
	\begin{equation} \label{STA} \tag{\textbf{STA}} E_0 (\eta) = E_0 (\psi\cup\mathcal{A}) \quad \text{for all  $\psi\subseteq\eta \cap \BB $ with $\psi\supseteq\theta  $ and all finite $\mathcal{A} \subseteq \R^d \setminus \mathcal{B} $}, 
\end{equation}
where $E_0(\cdot)$, the set of in- and out-neighbors of the origin, was defined in \eqref{equation_in_out_neighbors}.
 \item For every  $m \ge 1 $, there exists a subset of finite configurations consisting of precisely  $m $ elements $N_m\subseteq\mathbf N $ which is a zeroset with respect to the $dm $-dimensional Lebesgue measure and that has the property that for $\varphi\in\mathbf N\setminus N_m$ consisting of  $m $ elements, the set of out-neighbors  $\EE $ is  \emph{continuous}. Setting  $\NN\coloneqq \cup_{m \ge 1} N_m $, this means that for a finite $\varphi\in\mathbf N\setminus \NN $, there exists  $\delta>0 $ such that for every $x,y\in\varphi $ and every sequence  $ (z_w)_{w\in\varphi} \subseteq B_\delta (0) $,
	\begin{equation} \label{CON} \tag{\textbf{CON}} y+z_y\in \EE_{x+z_x} (\{w + z_w\co w\in\varphi\})\quad\text{if and only if} \quad y\in\EE_x (\varphi).
	\end{equation}  
		This assumption excludes finite configurations for which the graph is sensitive to small shifts of single or multiple nodes.   
\item There exists  $c_\ms{INF} > 0 $ with the following property:  let  $\psi\in\mathbf N_0 $ with  $\#\psi\ge c_\ms{INF}  $ and $\theta\in\mathbf N $.  We demand that
\begin{align} 
\label{INF} \tag{\textbf{INF}} \begin{split} &\EE (\psi)\subseteq\EE (\psi\cup\theta)\quad\text{if and only if}  \quad\EE (\psi)\subseteq\EE (\psi\cup\{y\} )\text{ for all } y\in\theta. \end{split}  
\end{align} 
In words, if the configuration  $\theta $ is such that no edges are removed by adding an element from  $\theta $ to $\varphi $, then also adding the entire set $\theta $ does not remove any edges  (and vice versa).   \end{enumerate} 

Each of these properties stays true if we increase $c_\ms{FIN},c_\ms{FIN2}  $ or  $c_\ms{STA}  $. Thus, we can set 
 $$
c_\ms{max} \coloneqq\max\{c_\ms{DEG},c_\ms{FIN}, c_\ms{FIN2},c_\ms{STA},c_\ms{INF} \} 
 $$
 and use it instead, where $c_\ms{DEG}$ represents a bound on the maximal node degree that is deduced in bullet point 4.\ below.  

 We now provide more detailed explanations for the conditions and their necessity.   
 \begin{enumerate}  
\item The scale invariance is a fundamental ingredient for controlling the asymptotic behavior of long edges.  This condition is satisfied by a variety of spatial networks such as the DT, the GG, and the RNG.   
\item[2./3.]  Condition \eqref{FIN}  is violated by the DT, the GG, and the RNG.  Moreover, if a graph does not fulfill condition  \eqref{FIN2}, then configurations may be possible with many points having very large edge lengths.  The RNG  (and therefore the DT and the GG) does not satisfy \eqref{FIN2}.  In this case, many nearby points with large combined edge lengths are possible by having two layers of points almost parallel to each other, as we elaborated in the paragraph after Equation \eqref{equation_representation_functional_X}. 
\item[4.]
In contrast to the stabilization conditions in \cite{hirsch,weakLLN}, we use a very specific class of stabilization regions $\BB$ based on cones. Nevertheless, it is still encompassed by more examples of spatial networks (such as RNG).  Our variant of the stabilization condition both allows for arbitrary modifications of the configuration outside the stabilization region $\BB$ but also adding points from the original configuration $\eta$ within $\BB$.
Our stabilization condition \eqref{STA} implies an alternative weaker version that is encompassed by even more examples of spatial networks (such as DT and GG), which is closer to the notion of stabilizing appearing in \cite{hirsch, weakLLN}. Namely, keeping the notation $\mc S_i$ from \eqref{STA} and all assumptions made there, we can define a \emph{stabilization radius} 
		 \begin{align}
			 \label{rr_eq}
\RR \colon \mathbf N_0\rightarrow [0,\infty],\ \varphi\mapsto  \max_{i\le I_d} \SS_i (\varphi), 
		 \end{align}
so that for all finite $\mc{A} \subseteq \R^d \setminus B_{\RR(X^*)}(0)$, setting $X^*=X\cup\{0\}$, we have
\begin{equation}\label{equation_weaker_stabilization_req}
E_0(X^*) = E_0((X^* \cap B_{\RR(X^*)}(0)) \cup \mc{A}).
\end{equation}
Defining stabilization by demanding the existence of an almost surely finite random variable, the stabilization radius $\RR(X^*)$, such that \eqref{equation_weaker_stabilization_req} is fulfilled is very similar to stabilization as it occurs in \cite{hirsch} and \cite{weakLLN}.

Additionally, \eqref{STA} yields a bound on the maximal node degree. In particular, when choosing $c_\ms{DEG} = I_d c_\ms{STA}$, we see that $\# E_0(\varphi)\le c_\ms{DEG}$ for all $\varphi\in\mathbf{N}_0$.  The implied uniformly bounded node degree helps to limit the number of edges that can contribute substantially to the power-weighted sum of edge lengths. 

The requirement that the lateral boundaries of the cones must  not be parallel to any of the axes is of technical nature and necessary in the proof of Lemma \ref{lemma_weak_law_of_large_numbers}. There we use a weak law of large numbers for Poisson functionals from \cite[Theorem 2.1]{weakLLN} which does not allow for points to be considered in the functional without their own scores contributing to the total sum. This can cause issues if we desire to compute probabilities which involve cones only containing a limited number of nodes up until a certain radius if the respective apex of the cone is close to the boundary of the observation window, as it happens in the proof of Lemma \ref{lemma_weak_law_of_large_numbers}. Here, we could imagine that there is some potential room to improve \eqref{STA} and drop the requirement about the lateral boundaries of the cones. For instance, one could try to be more lenient in a weak law of large numbers and also allow the consideration of points whose scores do not contribute. Another option would be to try to make use of the fact that \cite[Theorem 2.1]{weakLLN} allows for inhomogeneity of the points in some finer arguments. However, it is not clear if there are interesting examples of graphs that fulfill all other conditions but do not allow for lateral boundaries of the cones that are not parallel to any of the axes in \eqref{STA}.

\item[5.] In the theory of large deviations, it is common to make continuity assumptions in order to obtain asymptotically matching upper and lower bounds for the probability of rare events.  For instance, for the kNNs we want to avoid configurations where two distinct pairs of points have the same distance. 
 
\item[6.] $c_\ms{INF}  $ is necessary to ensure that the later introduced optimization problem that determines the rate of the large deviations, is indeed meaningful. In the simplest case, we would like to avoid situations in which there is a region for which adding a single point anywhere in it does not interfere with any existing edges but all of a sudden, a second point added to the region deletes one of the original edges. This could for example happen in the directed kNN with $k=2$ if the initial configuration consists of less than three points.  \end{enumerate}

\begin{figure} [h!]
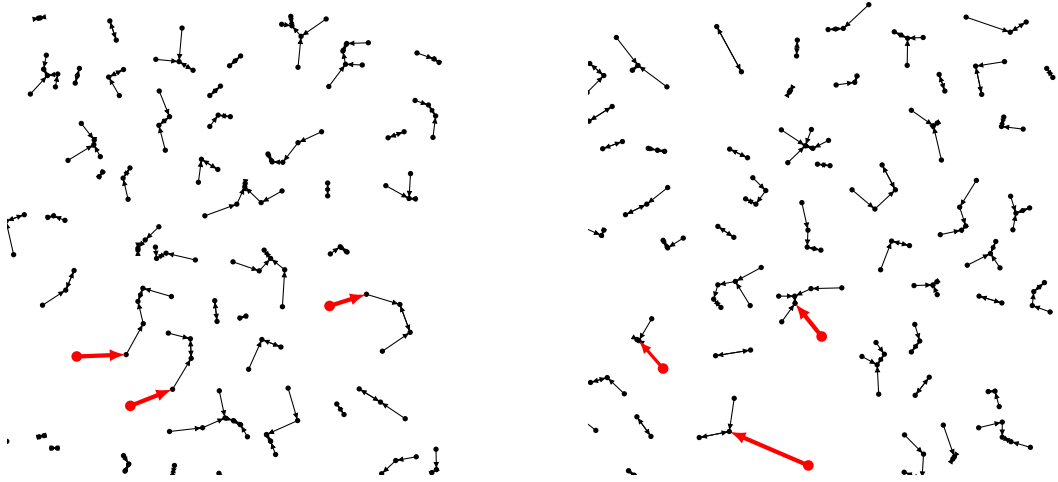
 \centering \input{Tikz_pictures/typicalH}  \hspace{1cm} \input{Tikz_pictures/largeH} 
\caption{Two configurations that result in a typical sum  (left) and an exceptionally large sum  (right)   of  $\a $-power weighted edge lengths with  $\a =15 $.  In each configuration, the three vertices inside an observation window with the most distant nearest neighbor are highlighted.  } \label{fig:simulations}  \end{figure}

%%LDP%
 Before introducing the deterministic optimization problem connected with the upper tails, Figure \ref{fig:simulations}  illustrates how the upper large deviations of  $H_n $ feature a condensate for the nearest neighbor graph.  There appears to be one large edge that carries the entire excess weight. 

The rate function in the large volume asymptotics will be given as a solution of an optimization problem.  To make this precise, we define the influence zone
\begin{equation}
A (\varphi,\psi)\coloneqq\big\{y\in\R^d \colon \EE_x (\psi)\not\subseteq \EE_x ( (\psi\cup\{y\} )) \text{ for some } x\in\varphi \cup\cup_{z\in\varphi}  \EE_z (\psi)\big\}
\end{equation}
 for configurations $\varphi\subseteq\psi\in\mathbf N $.  Loosely speaking, the cost of observing a certain configuration $\psi$ in the large-volume limit comes from the requirement that the influence zone may not contain any additional Poisson points. For instance, in the case of the kNN the influence zone describes the region of points, where adding an additional Poisson point would change one of the $k$ nearest neighbors of either an element of $\varphi$ or of a point that is itself one of the $k$ nearest neighbors of some element of $\vp$. 

 To be able to apply  \eqref{CON}  in Section \ref{subsec_low}, we set  $D (\varphi)\coloneqq\{y\in\R^d\colon\varphi\cup\{y\} \in N_{\#\varphi+1} \}  $ for a finite  $\varphi\in\mathbf N\setminus\NN $ as well as  $D _m' \coloneqq \{\psi\in\mathbf N\co\#\psi=m,|D (\psi)|>0\}  $ for  $m\in\N $.  Letting  $N_m'\coloneqq N_m\cup{D} _m' $ and $\NN'\coloneqq\cup_{m\ge1}  N_m' \supseteq \NN$, we then define the set of admissible configurations over which we optimize. These are configurations, whose total contributed power-weighted edge lengths exceed  $1 $, i.e.,
\begin{equation}\label{definition_sets_optimization}
B\coloneqq \Big\{ (\varphi,\psi) \co \varphi\subseteq\psi\in\mathbf N\setminus\NN',\, c_{\ms{INF}}  \le \#\psi< \ff, \,\sum_{x\in\varphi} \xi^{ (\a)}  (\psi - x) \geq 1\Big\}. 
\end{equation}
The most likely realizations in the large-deviation asymptotics are then the result of a delicate trade-off. We search for configurations that lead to a small influence zone $A$ but simultaneously exhibit edges that are long enough to be in the admissible set $B$.

 Now, we can state the main theorem, where $\mu_\a \coloneqq \E[\xi^{ (\a)}  (X\cup\{0\} )] $ denotes the expected edge length contribution of one vertex.   
 \bet[Upper large deviations]\label{theorem_main_upper_tails}  Let  $\a > d $ and  $r > 0 $.  Let the directed edge set  $\EE $ be scale invariant and satisfy  \eqref{FIN}, \eqref{FIN2},  \eqref{STA},    \eqref{CON}  and  \eqref{INF}.  Then,
\begin{equation}\label{eq_main_theorem}
\lim_{n \to \ff} \f1{n^{d^2/\a}} \log\P (H_n > \mu_\a + r) =-\inf_{ (\varphi,\psi)\in B}  |A (\varphi,\psi)| r^{d/\a}. 
\end{equation}
 \ent

The statement of Theorem \ref{theorem_main_upper_tails} indicates the necessity of a power larger than the dimension.  The usual speed for large deviations caused by homogenization in the situation of funtionals of this type of spatial random networks is $n^d$. If for $\a<d$, the equality in \eqref{eq_main_theorem} was still satisfied, we would have a faster speed than in a homogenization regime, which is not very reasonable and already gives a hint as to why our arguments require $\a>d$.

Next, we are going to assert that if the optimization problem has a strictly positive solution, then with high probability, only a negligible proportion of nodes is responsible for the entire excess when conditioned on the unlikely event.  In some cases, we can prove a sharper statement in the sense that only finitely many points carry the excess weight.  To make this precise, we introduce additional notation.  For configurations  $\varphi\subseteq\psi\in\mathbf N $, we will consider the order statistics of  $\xi^{ (\a)}  (\psi-x) ,\, x\in\varphi$.  That is, we let $Z^{ (i)}  (\varphi,\psi) $ denote the  $i $th largest element among $\{\xi^{ (\a)}  (\psi-x)\} _{x \in \varphi}  $.  In the case $\varphi = X\cap Q_n $ and  $\psi = X $, we abbreviate  $Z_n^{ (i)}  \coloneqq Z^{ (i)}  (X\cap Q_n,X) $ for  $i\ge 1 $. Besides that, recall the definition of the floor function $\lfloor t\rfloor \coloneqq \max\{m\in\Z\colon m\le t\}$ for $t\in\R$. In Theorem \ref{theorem_conditioned_convergence}, we add a further condition,  demanding that the volume of the influence zone does not become arbitrarily small even if using many nodes. 

\bet[Condensation conditioned on rare event]\label{theorem_conditioned_convergence}  Under the same conditions as in Theorem \ref{theorem_main_upper_tails}  and the additional assumption that $\inf_{ (\varphi,\psi)\in B}  |A (\varphi,\psi)| > 0 $, the following hold. \begin{enumerate}  
\item[a)] Let  $\eps\in (0, (1-d/\a)/ (2\a)) $ and  $\delta>0 $. Then, 
 $$
\P\bigg (\Big| (rn^d)^{-1} {\sum_{i\le \lfloor n^{d^2/\a-\eps} \rfloor}  Z_n^{ (i)}}  - 1\Big| > \delta\,\bigg|\, H_n>\mu_\alpha + r\bigg)\overset{n\to\infty} {\longrightarrow}  0. 
 $$
 
\item[b)]Additionally, assume there exists  $m_0 \ge 1 $ such that for every  $\delta\in (0,1) $,\begin{equation} \label{inequality_condensation_uniqueness} \inf_{ (\varphi,\psi)\in B} |A (\varphi,\psi)| <\inf_{ (\varphi,\psi)\in B,\, \sum_{i\le m_0}  Z^{ (i)}  (\varphi,\psi) < 1-\delta} |A (\varphi,\psi)|.   \end{equation}  Then, for every  $\delta>0 $, 
 $$
\P\bigg (\Big| (rn^d)^{-1} {\sum_{i\le m_0}  Z_n^{ (i)}}  -1\Big| > \delta\, \bigg|\, H_n>\mu_\alpha + r\bigg)\overset{n\to\infty} {\longrightarrow}  0. 
 $$
\end{enumerate}  \ent 
Condition \eqref{inequality_condensation_uniqueness}  implies that any optimal configuration consists of at most $m_0$ nodes. As will be shown in Section \ref{subsec_NNG}, the nearest neighbor graph (NNG) for large  $\a $ is an example for a graph satisfying \eqref{inequality_condensation_uniqueness}.

\begin{remark}  Theorem \ref{theorem_main_upper_tails}  also can be applied if instead of  $X $, we consider a Poisson process  $Y $ with intensity $n^{-\beta d}  $ for some  $\beta < 1 $.  Scaling  $Y $ by  $n^{-\beta}  $ yields a Poisson process with intensity 1 and the window becomes  $Q_{n^{1 -\b}}  $.  The mean is given as  $\E\big[\xi^{ (\a)}  \big( (n^{-\beta} Y)\cup\{0\} \big)\big] =n^{-\a\beta} \mu_\a, $ finally yielding upper tails of the form
	\begin{align*}  &\lim_{n \to\ff} \f1{n^{ (1-\beta)d^2/\a}} \log\P\bigg (\f1{n^{d+\beta (\alpha-d)}} \sum_{x\in Y\cap Q_n}  \xi^{ (\a)}  (Y-x) > \mu_\a + r\bigg) \\&=\lim_{n \to\ff} \f1{n^{ (1-\beta)d^2/\a}} \log\P\bigg (\f1{n^{d-\beta d}} \sum_{x\in X\cap Q_{n^{1- \b}} }  \xi^{ (\a)}  (X-x) > \mu_\a +r\bigg) = -\inf_{ (\varphi,\psi)\in B}  |A (\varphi,\psi)|r^{d/\a}, 
	\end{align*}  
where in the last line we applied Theorem \ref{theorem_main_upper_tails}  with  $n'=n^{1-\beta}  $.   \end{remark} 

\begin{remark}
Another interesting graph to examine in terms of a condensation phenomenon is the directed spanning forest (DSF). Very loosely speaking, this graph draws an edge from a node to the closest other node that has a higher value in the $d$th coordinate, see \cite{spanningForest}. This graph does not satisfy condition \eqref{FIN} required for the upper large deviations and condensation. Further, in the given form of the DSF, this would be one of the few common examples where \eqref{STA} is violated due to the lateral boundary part. Nevertheless, we suspect the total power-weighted edge lengths for $\a>d$ for the DSF to admit upper large deviations with a condensate that might even involve the same optimization problem as it appears in Theorem \ref{theorem_main_upper_tails}. One would need a more generous concentration bound that does not rely on \eqref{FIN} to proof Lemma \ref{lemma_hoeffding} and as pointed out in the explanation of \eqref{STA}, we are also confident that it is possible, with finer arguments, to drop the lateral boundary condition from \eqref{STA}. Here, this issue could even be avoided if the search process of the DSF for the closest point would not be parallel to one of the axes.
\end{remark}

\begin{remark}
We limit ourselves to the study of the functional representing power-weighted edge lengths of spatial random networks in terms of its upper large deviations. Even the consideration of this functional for a power larger than the dimension restricts the class of admissible graphs heavily. Nevertheless, we can imagine that there is room to potentially improve this and, on top of the graph, generalize the functional as well. An idea would be to consider functional-graph combinations that for a node to have a large score would require a relatively large region to contain no or only a limited amount of points. This would include the total sum of power-weighted edge lengths for the kNN and $\beta$-skeleton. An example besides our studied functional that would fit this description could be the sum of power-weighted circumradii of the simplices in the DT. However, if we, like in this specific example, study condensation phenomena for functionals that we apply to the DT, we would run into other issues that were described in the explanations of our conditions.
\end{remark}
\section{Applications of Theorem \ref{theorem_main_upper_tails}}\label{sec_application} We  verify that the (un-/bidirected) kNN and suitable  $\b $-skeletons satisfy the conditions in Theorem \ref{theorem_main_upper_tails}. 

\subsection{$k$-nearest neighbor graphs}\label{subsec_kNN} 
In the \emph{kNN}, a directed edge is drawn from each node to the  $k\ge1 $ points that are closest in Euclidean distance. As explained in Section \ref{sec_model}, this definition gives rise to undirected and bidirected kNNs. For  $j\le k $, we define the distance from the origin to the  $j $th closest point in a configuration by   
$$
\mc D_j \colon \mb N_0 \rightarrow [0,\infty),\ \varphi \mapsto \inf\{r>0\co \varphi(B_r(0)) \ge j+1\}. 
$$
 This leads to the set of the  $k $ nearest neighbors of the origin  
\begin{equation}\label{definition_edges_kNN}
\mc E\colon \mb N_0 \rightarrow \mb N,\ \varphi \mapsto \{x\in\varphi\cap B_{\mc D_k(\varphi)}(0)\}\setminus\{0\}.
\end{equation}
 We will use the lexicographical order to determine the  $k $ nearest neighbors of a node in case more than  $k $ neighbors are potential candidates. 
 In the following, we quickly verify the conditions in Theorem \ref{theorem_main_upper_tails}.  \begin{enumerate} 
\item $\EE$ defined as in \eqref{definition_edges_kNN} is scale invariant.
\item The bounded node degree \cite[Lemma 8.4]{yukichProbabilityBook}  entails \eqref{FIN} with  $c_\ms{FIN}=c_\ms{DEG} $ since all nodes that are affected by adding a new vertex to the configuration must be part of an edge with the new vertex. 
\item Let $\varphi\in\mathbf{N}$ and $M>0$ be arbitrary. To ease presentation, we consider $k=1$ first. Each vertex $x\in\varphi\cap B_M(0)$ incident to an edge longer than $M$ defines a ball of radius at least $M$, centered at $x$, that does not contain any other vertices in its interior. Hence, scaling the radii by 1/2 gives rise to a family of balls that are pairwise disjoint, each having radius at least $M/2$. Thus, the number of nodes within $B_M(0)$ that are incident to an edge larger than $M$ is at most $|B_{2M}(0)|/|B_{M/2}(0)|=4^d$. 

	Now, let $k \ge 2$ be general and set $\varphi'=\varphi$. Starting with a node
	$$x\in\argmax_{z\in\varphi'\cap B_M(0)} \{|y|\co y\in\EE(\varphi-z)\text{ and }|y|>M\},$$
	we delete all points in $\varphi'$ that are within the interior of $B_{\mc D_k(\varphi-x)}(x)\setminus\{x\}$, which are at most $k-1$, and mark $x$ as already dealt with. We repeat this procedure recursively, ignoring nodes in the index of the $\argmax$ that are already marked, until all nodes in $\varphi'\cap B_M(0)$ are either marked or not associated with an edge of length exceeding $M$. Then, by the same arguments as in the case $k = 1$, the interiors of the balls $B_{\mc D_1(\varphi'-x)/2}(x)$ are pairwise disjoint for $x\in\{z\in\varphi'\cap B_M(0)\co |y|>M \text{ for some }y\in\EE(\varphi-z)\}$ and $\#\{z\in\varphi'\cap B_M(0)\co |y|>M \text{ for some }y\in\EE(\varphi-z)\}$ is bounded by $4^d$. Moreover, for every marked node left in the thinned configuration $\varphi'\cap B_M(0)$, we deleted at most $k-1$ nodes from $\varphi$ and thus, we deduce that the total number of nodes in $\varphi\cap B_M(0)$ incident to an edge of length exceeding $M$ is at most $c_\ms{FIN2}\coloneqq k4^d$ which yields \eqref{FIN2}.
\item Considering only the undirected kNN, from \cite[Lemma 6.1]{CLT} it follows that we can find a collection of cones such that  $\RR $ can be used as stabilization radius in the weaker sense of \eqref{STA} with $c_\ms{STA} \coloneqq k+1 $. Now, for $\varphi\in\mathbf{N}_0$, let $P_i$ denote the set of the $c_\ms{STA}$ closest points to the origin in $\varphi\cap (S_i\setminus\{0\})$. If the intersection does not contain $c_\ms{STA}$ points, then let $P_i=\varphi\cap S_i\setminus\{0\}$ or if there are more than $c_\ms{STA}$ candidates let the lexicographical order decide which of the candidates the furthest away to include in $P_i$ and put $\theta\coloneqq\{x\in P_i\co i\in\{1,\dots,I_d\}\}$. Then,  \cite[Lemma 8.4]{yukichProbabilityBook}, which asserts that the undirected kNN has bounded node degree, and its proof imply that we can choose the cones in a way such that for this choice of $\theta$, the condition \eqref{STA} is satisfied. Further, because \eqref{STA} only incorporates $\EE$, condition \eqref{STA} follows for the undirected, bidirected and directed kNN.
\item The continuity condition \eqref{CON} is satisfied with
$
N_m \coloneqq \{\varphi\in\mb N\co \#\varphi=m \text{ and }|w-x|=|y-z|>0 \text{ for some } w,x,y,z\in\varphi \text{ with }\{w,x\}\neq\{y,z\}\} 
$
 as the set of configurations containing  $m $ nodes, where there are pairs of nodes with equal distances. 
\item We choose  $c_\ms{INF}\coloneqq k+1 $ to ensure that each node has  $k $ neighbors. Then, \eqref{INF} is satisfied since for  $\varphi\in\mb N_0 $ with  $\#\varphi\ge c_\ms{INF} $, a node in the set  $\mc E(\varphi) $ only vanishes when adding a vertex within the interior of the ball  $B_{\mc D_k(\varphi)}(0) $. Adding more vertices can only cause more differences. \end{enumerate}

	%
	%BETA SKEL
	%
 \subsection{{ $\b $}-skeleton}  
\begin{figure}[h!] \centering 
	 \begin{tikzpicture}[scale=1.2] \coordinate (x) at (0,0); \coordinate (y) at (3,0); 

 \draw[color=red] (y) arc[start angle={-18.43}, end angle={180+18.43},radius={sqrt(1.5^2+0.5^2)}]; \draw[color=red] (x) arc[start angle={180-18.43}, end angle={360+18.43},radius={sqrt(1.5^2+0.5^2)}]; 

 \coordinate (z) at (1,2); 

 \draw[dashed,color=blue] (x) -- (z); \draw[dashed,color=blue] (y) -- (z); \path (x) -- (z) -- (y) pic["\tiny{ $\arcsin\f1{\b} $}",draw=blue,{Latex[scale=0.8]}-{Latex[scale=0.8]},angle eccentricity=1.5,angle radius=0.5cm] {angle=x--z--y}; 

 \fill (x) circle (1pt) node[xshift=-0.25cm] { $x $}; \fill (y) circle (1pt) node[xshift=0.25cm] { $y $}; \draw (x) -- (y); \end{tikzpicture} \hspace*{2cm} 
	 \input{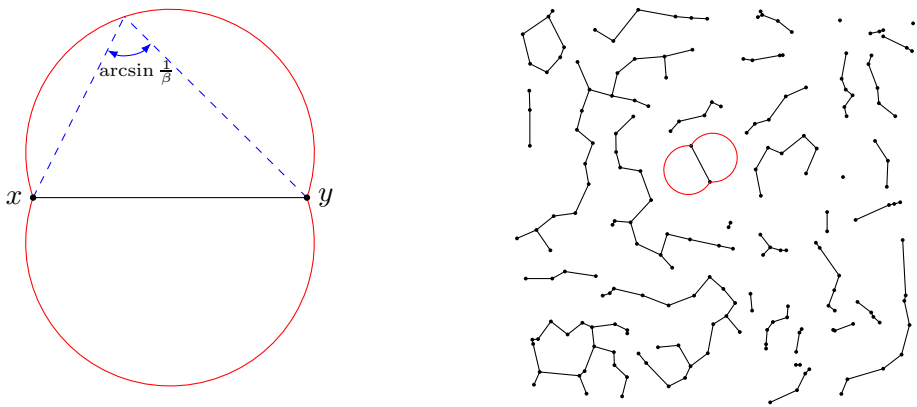} \caption{Illustration of an edge in the  $\b $-skeleton and a random simulation of the  $\b $-skeleton with  $\b=1.2$.}\label{fig:b_skeleton} 
 \end{figure}

\emph{ $\b $-skeletons} are geometric graphs that are popular in applications in pattern recognition \cite{kirkpatrick} and machine learning \cite{toussaint}.  The 2D  $\b $-skeleton, $\b > 1$, has an edge between two nodes  $x $ and  $y $ if there is no vertex that has an angle, generated by the two lines to  $x $ and  $y $, that is larger than  $\gamma\coloneqq\arcsin(\b^{-1}) $. In other words, there is an edge if the union $C(x,y)$ of the two disks with radius  $\b|x-y|/2 $ and having   $x $ and  $y $ on their boundary does not contain any other vertices, see Figure \ref{fig:b_skeleton}.  This construction rule determines the set of neighbors $\mc E $. Note that this definition also makes sense for $\b = 1$, leading to a spatial network known as \emph{Gabriel graph}.

Although the $\b$-skeleton can also be defined in higher dimensions, we henceforth restrict our attention to the 2D  $\b $-skeleton for two reasons. First, the 2D case already covers the vast majority of applications of the  $\b $-skeletons. Second, as we will see below, already in the 2D case, the verification of condition \eqref{FIN} requires delicate geometric arguments. Although we believe an extension to higher dimension is possible, this would entail an even more tedious geometric analysis. Since the focus of our article is on presenting novel probabilistic aspects of large deviations in a geometric context, it would not be appropriate to devote several pages of trigonometry arguments to the verification of the conditions in three and higher dimensions.

We now verify that the $\b$-skeletons satisfy the conditions of Theorem \ref{theorem_main_upper_tails}. To that end, we state an auxiliary result capturing the stabilization properties of $\b$-skeletons needed for condition \eqref{STA}. Since the $\b$-skeleton is intrinsically an undirected graph, we henceforth consider all appearing edges as undirected in order to make the presentation more accessible.
	 \bel[Stabilization for $\b$-skeletons]\label{lemma_b_skel_cones} For  $\b> 1 $, there is a collection of cones  $(S_i)_{1\le i\le I_2} $ satisfying the requirements of \eqref{STA} with $c_\ms{STA} = 2$. 
	 \enl 
	 \bep
	 We choose the cones $S_i$, $i \le I_2 $ sufficiently thin and not axes-parallel such that for any  $r>0 $, the angle generated by starting from the origin, proceeding to any point in  $S_i\cap B_r(0) $ and ending at any point in  $S_i\cap \partial B_r(0)$ exceeds $\gamma $. Now, if $x \in \vp$ is the closest point to $0$ contained in $S_i$, then $\angle 0xy  > \g$ for every $y \in S_i$ with $|y| \ge |x|$ and $x\neq y$. Thus, there cannot be an edge between the origin and  $y $.
	 
	 To construct $\th$, we first let $P_i$ denote the closest point to the origin in $\vp \cap (S_i\setminus\{0\})$ if the intersection is non-empty (resolving potential ties by choosing the lexicographic minimum). Then, we put $\th \coloneqq \{P_i \co \vp \cap S_i \ne\es\}\cup\{0\}$.
	 \enp 
Leveraging Lemma \ref{lemma_b_skel_cones}, we now verify conditions 1, and 4--6. The application of Theorem \ref{theorem_conditioned_convergence} for the  $\b $-skeleton is verified in Section \ref{sec_application_condensation} below.

\been
\item  $\mc E $ for the $\b $-skeleton, where  $\b>1 $, is scale invariant. 
\item[4.] This is the content of Lemma \ref{lemma_b_skel_cones}.
\item[5.] The continuity condition \eqref{CON} is satisfied with
$
N_m \coloneqq \{\varphi\in\mb N\co \#\varphi=m \text{ and }\vp\cap \partial C(x, y) \ne \{x, y\}\text{ for some } x,y\in\varphi \} 
$
 as the set of configurations containing  $m $ nodes, where there are two nodes that have a vertex on the boundary of the union of disks illustrated in Figure \ref{fig:b_skeleton}. 
\item[6.] To remove a $\b $-skeleton edge $e$, only one node in $C(e)$ is sufficient. Hence, $c_\ms{INF}=1 $.  
\enen 

In the rest of this section, we verify conditions \eqref{FIN} and \eqref{FIN2}.

 For  $e_1, e_2\in\R^2 $ and the edge $e := (e_1, e_2)$ with  $|e_1-e_2|\ge a>0 $, we define the point between  $e_1 $ and  $e_2 $ that has distance  $a $ from  $e_1 $ by  $h_a(e) \coloneqq e_1+(e_2-e_1)a/|e| $. Further, let  $M(e) $ be a point at distance  $\b|e|/2 $ to both,  $e_1 $ and  $e_2 $. In other words,  $M(e) $ represents the center of one of the two disks that unioned represent  $C(e)$, see Figure \ref{fig:b_skeleton}. In some cases, we will need to make a specific choice between one of the two options, and then we will state this clearly. Finally, let $\Delta_{M(e)}(e) $ be the triangle formed by $M(e) $ and  $e $.

 \bel[Disjoint regions for  $\b $-skeletons]\label{lemma_b_skel_disjoint_disks} 
 Let $e_1, e_2, f_1, f_2 \in \R^2$ be pairwise distinct, and assume that   $e = \{e_1, e_2\}, f = \{f_1, f_2\} \in E(\{e_1, e_2, f_1, f_2\}) $. Then, 
 \begin{enumerate} 
\item[i)]  $f $ does not intersect  $\Delta_{M(e)}(e)$;
\item[ii)] there exists a constant  $c_\ms{disj} = c_\ms{disj}(\b) \in(0,1/2) $,  such that if  $|e| \wedge |f| \ge a $ for some $a>0 $, then 
$$
B_{c_\ms{disj}a}(h_m(e)) \cap B_{c_\ms{disj}a}(h_{m'}(f)) = \emptyset 
$$
 for all  $m\in[a/2,|e|-a/2] $ and  $m'\in[a/2,|f|-a/2] $.  
 \end{enumerate} \enl  

 We postpone the proof of Lemma \ref{lemma_b_skel_disjoint_disks} to the end of this section, and elucidate how to verify condition \eqref{FIN2}. First, instead of bounding the number of nodes in $B_M(0)$ incident to a long edge, we may bound the number of disjoint long edges with one endpoint in $B_M(0)$. Then, we apply Lemma $\ref{lemma_b_skel_disjoint_disks}$ for every pair of such disjoint edges $e$ and $f$ with $a \coloneqq M$, $m \in\{a/2,|e|-a/2\}$ and $m'\in\{a/2,|f|-a/2\}$, depending for which choice of $m$ and $m'$ the points $h_m(e)$ and $h_{m'}(f)$ are closer to $B_M(0)$. Hence, having $k \ge 1$ disjoint long edges with an endpoint in $B_M(0)$ leads to $k$ disjoint disks with radius $c_\ms{disj}M$ that are contained entirely within $B_{2M}(0)$. Thus, the number of such edges is at most $|B_{2M}(0)|/|B_{c_\ms{disj}M}(0)| = |B_2(0)|/|B_{c_\ms{disj}}(0)|$.

 \begin{wrapfigure}{r}{0.405\textwidth}
 \centering \begin{tikzpicture}[scale=1.1] \coordinate (x) at (0,0); \coordinate (y) at (5,0); \coordinate (z) at (2.5,-0.5); \coordinate (w) at (2.3,-2.985); 

 \coordinate (a) at (1.75,-2.941); \coordinate (b) at (3.25,-2.941); \coordinate (e) at (2.5,-3.13); 

 \node[xshift=-0.2cm] at (x) { $\small e_1 $}; \node[xshift=0.2cm] at (y) { $\small e_2 $}; \node[yshift=-0.2cm] at (z) { $\small M(e) $}; 

 \node[color=red, xshift=-0.2cm] at (a) { $\small f_1 $}; \node[color=red, xshift=0.2cm] at (b) { $\small f_2 $}; \node[color=red, xshift=0.2cm, yshift=-0.2cm] at (e) { $\small M(f) $};

 \tkzDrawArc(z,x)(y) 

 \tkzDrawArc[red!50](e,a)(b)

 \coordinate (f) at (1.9,0); \coordinate (g) at (1.9,-0.375); \tkzDrawArc[blue](x,g)(f) \node[color=blue, xshift=-0.1cm, yshift=0.225cm] at (g) { $\small \t $}; 

 \draw[dashed] (z) -- (x); \draw[dashed] (z) -- (y); 

 \draw[dashed, color=red] (a) -- (e); \draw[dashed, color=red] (b) -- (e); 

 \draw[dotted, color=blue] (w) -- (x); \draw[dotted, color=blue] (w) -- (y); 

 \fill (z) circle (0.6pt); \fill[red] (e) circle (0.6pt); \fill[blue] (w) circle (1.0pt) node[xshift=-0.1cm, yshift=-0.2cm] { $\small y $}; 

 \draw[red] (a) -- node[above, yshift=-0.09cm] { $f $} (b); \draw (x) -- node[above, yshift=-0.09cm] { $e $} (y); 

 \end{tikzpicture} 
	 \caption{Illustration of the statement of Lemma \ref{lemma_b_skel_FIN} including the inserted node relevant for \eqref{FIN}. The extended line between  $f_1 $ and  $M(f) $ is tangent to the disk segment.}\label{fig:b_skeleton_FIN} 
 \end{wrapfigure}
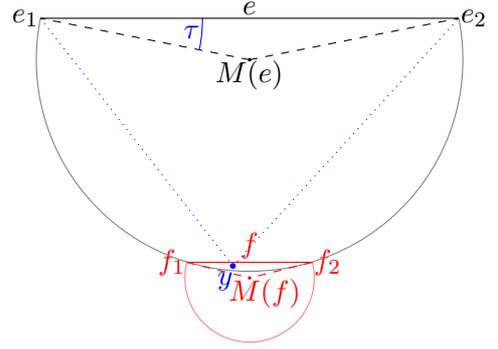

Finally, we verify condition \eqref{FIN}. To achieve this goal, note that the amount of edges that can arise from  $y \in \R^2$ is limited by the bound on the node degree. Hence, it remains to consider the number of edges removed by adding the point $y$. In particular, the number of disjoint edges removed is sufficient. Here, a key observation is that if $e$, $f$ are disjoint edges with $y \in C(e) \cap C(f)$, then this implies a very particular relative configuration for $e$ and $f$. More precisely, the edges $e$ and $f$ do not intersect, and the triangle $\De_y(e)$ does not contain an endpoint of $f$ and vice versa. Hence, if we consider the cones $S_y(e)$ and $S_y(f)$ with apex $y$ obtained by extending these triangles, then there are only 3 options: (i) $S_y(e) \cap S_y(f) = \{y\}$, (ii) $S_y(e) \su S_y(f)$, or (iii) $S_y(f) \su S_y(e)$. In the latter cases, we say that $e$ and $f$ are related. Since the angle at the apex of each of these cones is at least $\gamma$, the number of equivalence classes of related edges is at most $2\pi /\g$.

Hence, to complete the proof of condition \eqref{FIN} it suffices to bound the number of elements in each equivalence class.
 For this  step, we need two further results. To state them, we set $\t \coloneqq \arccos(\b^{-1})$.

 \bel[Exclusion of short edges]\label{lemma_b_skel_FIN} 
 Let $e_1, e_2, f_1, f_2 \in \R^2$ be pairwise distinct, and assume that   $e \coloneqq \{e_1, e_2\}, f \coloneqq \{f_1, f_2\} \in E(\{e_1, e_2, f_1, f_2\})$ and that $|f| \le \tan(\t) |e|$. Furthermore, let $y \in C(e)$ be such that  $f $ crosses  $\Delta_y(e) $ between  $e $ and  $y $. Then,  $y \in \Delta_{M(f)}(f). $
 \enl 
 The configuration in Lemma \ref{lemma_b_skel_FIN} is sketched in Figure \ref{fig:b_skeleton_FIN}. 
Next, for  $\varphi\in\mb N $,  $y\in\R^2 $ and $e \in E(\varphi)$ with $y \in C(e)$, we define 
 \begin{equation} 
	 E_\ms{REC}(\varphi,y, e) \coloneqq \{f\in E(\varphi) \co S_y(e) \su S_y(f) \text{ and }y \in C(f)\}
 \end{equation} 
 as the set of recorded edges.

 \bel[Size bound for recorded set]
 \label{lemma_b_skel_FIN2} 
 There exists  $c_\ms{edges} = c_\ms{edges}(\b)>0 $  such that for any  $\varphi\in\mb N$, $e \in E(\vp)$ and $y\in\R^2 $  with $y \in C(e)$, we have $\#E_\ms{REC}(\varphi,y, e) \le c_\ms{edges}$.
 \enl 

 Note that once Lemma \ref{lemma_b_skel_FIN2} is established, condition \eqref{FIN} is verified since then the total number of deleted edges is at most $c_\ms{edges}2\pi/\gamma$.
Hence, it remains to prove the auxiliary results Lemmas \ref{lemma_b_skel_disjoint_disks}, \ref{lemma_b_skel_FIN}, and \ref{lemma_b_skel_FIN2}. 

 \bep[Proof of Lemma \ref{lemma_b_skel_disjoint_disks}]\phantom{a}\\[-2ex]

{\bf Part i)} In the setting of Lemma \ref{lemma_b_skel_disjoint_disks}, assume that  $f $ intersects  $\Delta_{M(e)}(e) $ and note that the nodes  $f_1 $ and  $f_2 $ have to be outside $C(e)$ for  $e $ to exist. But since  $f $ intersects  $\Delta_{M(e)}(e) $, at least one of  $e_1 $ and  $e_2 $ is in  $B_{|f|/2}(h_{|f|/2}(f)) $. Therefore, $f $ would not exist in the GG, and thus also not in the $\b$-skeleton. Hence, $\Delta_{M(e)}(e) $ cannot intersect $f $.

{\bf Part ii)}
 Repeating the above argument for the second choice of $M(e)$ yields a rhombus with centroid  $h_{|e|/2}(e) $ that cannot be intersected by other edges. However, since the side lengths of this rhombus are of order  $|e|>a $, there exists a constant $c_\ms{disj} = c_\ms{disj}(\b)\in(0,1/2) $ such that any disk with center between  $h_{a/2}(e) $  and  $h_{|e|-a/2}(e)$ and radius $c_\ms{disj}a $ also has distance of more than  $c_\ms{disj}a $  to the boundary of the rhombus (and similarly for $e$ replaced by $f$). Since the rhombus linked to any edge cannot be intersected by another edge, it follows that the disk associated with  $e $ and the disk associated with  $f $ are disjoint.
 \enp

 \bep[Proof of Lemma \ref{lemma_b_skel_FIN}] Since $e$ is an edge in the $\b$-skeleton, the nodes   $f_1 , f_2 $ lie outside the interior of $C(e)$. We first consider the case where $f_1, f_2$ are contained in the boundary of $C(e)$, and the segments $[M(f), f_1]$, $[M(f), f_2]$ are tangent to $C(e)$. Then, $\Delta_{M(f)}(f)\cap C(e) $ yields a full circular segment of  $B_{\b|e|/2}(M(e))$ so that $y \in \Delta_{M(f)}(f)$. We assert that if $[M(f), f_1]$ and $[M(f), f_2]$ are tangent to $C(e)$, then  $|f|=\tan(\t)|e|$. Since $y \in \Delta_{M(f)}(f)$ will remain true if we shorten $|f|$, this will conclude the proof of the lemma.

 To prove that $|f| = \tan(\t)|e|$, note that the tangency implies that $M(e)f_1M(f)$ is a right triangle. Thus, 
$$
 \frac{|f_1-M(f)|}{|f_1 - M(e)|}= \frac{|f_1-M(f)|}{\b|e|/2} = \tan(\t). 
$$
 Next, also $f_1M(f) h_{|f|/2}(f) $ defines a right triangle so that  
$
{|f|}/{(2|f_1-M(f)|)} = \cos(\t) = \b^{-1}.
$
Finally, combining these two relations yields the asserted
$
|f| = \tan(\t)|e|. 
$
 \enp 

 \bep[Proof of Lemma \ref{lemma_b_skel_FIN2}] First, we note that $E_\ms{REC}(\varphi,y, e)$ contains at most one edge that is shorter than $\tan(\t)|e|$. Indeed, suppose that $f \ne f'$ are two such edges with $S_y(f) \su S_y(f')$.  Now, from Lemma \ref{lemma_b_skel_disjoint_disks} i), we know that $f'$ cannot intersect  $\Delta_{M(f)}(f) $ and therefore also not  $\Delta_y(e)\cap\Delta_y(f) $. This contradicts Lemma \ref{lemma_b_skel_FIN}.

 Hence, it suffices to bound the number of $f \in E_\ms{REC}(\varphi,y, e)$ with $|f| \ge \tan(\t)|e|$. To achieve this goal, let $f^{(1)},\dots,f^{(c)}\in E_\ms{REC}(\varphi,y,e)$ be disjoint edges, each of length at least $\tan(\t)|e|$. Note that none of these edges can intersect. Further, for the edge  $e $ to exist, the edges  $f^{(1)},\dots,f^{(c)} $ must also fully cross the disk segment  $C(e) $ as it is drawn in Figure \ref{fig:b_skeleton_FIN}.

Then, for all $i \le c$, the edge $f^{(i)}$ crosses the cone $S_y(e)$ somewhere since $S_y(e)\subseteq S_y(f^{(i)})$. In particular, $f^{(i)}$ has to cross the triangle $\De_y(e)$. If that was not the case and $f^{(i)}$ would cross $S_y(e)\setminus \De_y(e)$, then $e_1,e_2\in \De_y(f^{(i)})\subseteq C(f)$ which would contradict the existence of the edge $f^{(i)}$. It is impossible for $f^{(i)}$ to cross both $S_y(e)\setminus \De_y(e)$ and $\De_y(e)$ because then it would have to intersect $e$. 

 Then, by Lemma \ref{lemma_b_skel_disjoint_disks} each $f^{(i)}$ generates a disk with radius  $c_\ms{disj}\tan(\t)|e| $ with center that has to be within distance  $\tan(\t)|e| $ of  $C(e) $, disjoint from the disks created by other edges larger than  $\tan(\t)|e| $. Thus, the total number of long edges that can cross  $\Delta_y(e) $ is bounded by 
\begin{equation*}
	\label{equality_b_skel_edges_triangle} \frac{2|B_{2\tan(\t)|e|+\b|e|/2}(M(e))|}{\pi(c_\ms{disj}\tan(\t)|e|)^2}=\frac{2|B_{2\tan(\t)+\b/2}(0)|}{\pi c_\ms{disj}^2\tan(\t)^2} \eqqcolon c_\ms{edges}(\beta)-1,
 \end{equation*}
 thereby concluding the proof.

	\enp

\section{Applications of Theorem \ref{theorem_conditioned_convergence} a) and b)}\label{sec_application_condensation}
In this section, we verify the conditions of Theorem \ref{theorem_conditioned_convergence} a) for the graphs from Section \ref{sec_application}. We also apply Theorem \ref{theorem_conditioned_convergence} b) to the NNG. To ease the overall presentation, we start with the latter.

\subsection{Theorem \ref{theorem_conditioned_convergence} b) for the NNG}\label{subsec_NNG}

We start with an auxiliary result simplifying the definition of the influence zone for the NNG. Loosely speaking, we can ignore the constraints on the outneighbors of $\vp$ and can concentrate on the areas influencing the nearest neighbors of points in $\vp$ itself. 
\bel[Influence zone for the NNG]\label{lemma_optimization_problem_for_NNG}
It holds that
$$\inf_{(\varphi,\psi)\in B} |A(\varphi,\psi)| = \inf_{(\varphi,\psi)\in B} |\cup_{x\in\varphi} B_{\mc D_1(\psi-x)} (x)|.$$
\enl

\begin{remark}\label{remark_optimization_problem_for_NNG}
	An adaptation of the proof of Lemma \ref{lemma_optimization_problem_for_NNG} shows that it remains true if on both sides we replace $B$ by $\{(\varphi,\psi)\in B\co\sum_{i\le m_0} Z^{(i)}(\varphi,\psi) < 1-\delta\}$. The proof can be replicated without significant alterations.
\end{remark}

Next, we further examine the geometric interpretation of the optimization problem.
\bel[One single large ball is the unique optimal solution for the NNG and $\a\gg d$]\label{lemma_one_ball_asymptotically_optimal}
There exists $\a_0 > d$ such that the configuration $(\{{0}\},\{{0},(1,0\dots,0)\})$ solve the optimization problem for all $\a\ge\a_0$. In particular,
\begin{equation}\label{equality_one_ball_asymptotically_optimal}
\inf_{(\varphi,\psi)\in B} |A(\varphi,\psi)| = |B_1(0)| = \k_d.
\end{equation}
Moreover, for every $\delta>0$ there exists $\eps>0$ such that 
$|\cup_{x\in\varphi} B_{\mc D_1(\psi-x)}(x)| \ge (1+\eps)\k_d$
holds for all $(\varphi,\psi)\in B$ with $\max_{x\in\varphi}\mc D^{(\a)}_1(\psi-x) < 1-\delta$.
\enl

Hence, to verify the application of Theorem \ref{theorem_conditioned_convergence} part b) for the NNG, only the proofs of Lemmas \ref{lemma_optimization_problem_for_NNG} and \ref{lemma_one_ball_asymptotically_optimal} are necessary.

\bep[Proof of Lemma \ref{lemma_optimization_problem_for_NNG}]
First, by the definition of $A$ in the case of the NNG, we have that
$$|A(\varphi,\psi)| = |\underbrace{\cup_{x\in\varphi} \big(B_{\mc D_1(\psi-x)} (x)\cup\cup_{z\in\mc{E}_x(\psi)} B_{\mc{D}_1(\psi-z)}(z)\big)}_{\eqqcolon K(\varphi,\psi)}|$$
for all $(\varphi,\psi)\in B$, since in the NNG an edge can only be deleted if an additional node is put within the open ball with radius given by $\mc{D}_1(\cdot)$, centered at a vertex in $\psi$. This implies $\inf_{(\varphi,\psi)\in B} |A(\varphi,\psi)| \ge \inf_{(\varphi,\psi)\in B} |\cup_{x\in\varphi} B_{\mc D_1(\psi-x)} (x)|$. 

For the other direction, let $\eps>0$ and $(\varphi,\psi)\in B$ be arbitrary.
Now for $\delta>0$ we introduce an extended configuration $\theta_{\delta}\supseteq\psi$ by adding a further point to $B_\delta(x)\setminus\{x\}$ for all $x\in\cup_{z\in\varphi} (\mc{E}_z(X))\setminus\varphi$. Hence,
$$\underbrace{|\cup_{x\in\varphi} B_{\mc D_1(\psi-x)-\delta} (x)|}_{\overset{\delta\downarrow 0}{\longrightarrow} |\cup_{x\in\varphi} B_{\mc{D}_1(\psi-x)} (x)|} \le |K(\varphi,\theta_{\delta})| \le \underbrace{|\cup_{x\in\varphi} \big(B_{\mc{D}_1(\psi-x)} (x)\cup\cup_{z\in\mc{E}_x(\theta_{\delta})} B_{\delta}(z)\big)|}_{\overset{\delta\downarrow 0}{\longrightarrow} |\cup_{x\in\varphi} B_{\mc{D}_1(\psi-x)} (x)|},$$
where the convergences follow because the chosen configurations are finite. Thus, we can choose $\delta$ small enough for
$\big||K(\varphi,\theta_{\delta})|-|\cup_{x\in\varphi} B_{\mc D_1(\psi-x)} (x)|\big| \le \eps$. Scaling all the configurations with $1+\eps$ gives that $\sum_{x\in(1+\eps)\varphi} \xi^{(\a)}((1+\eps)\psi-x) \ge 1+\eps$. Note that due to the finiteness of the configurations in $B$, we can let $\delta$ be small enough such that still $\sum_{x\in(1+\eps)\varphi} \xi^{(\a)}((1+\eps)\theta_{\delta}-x) \ge 1$ which implies that $((1+\eps)\varphi,(1+\eps)\theta_{\delta})\in B$. Thus,
$$
	|\cup_{x\in\varphi} B_{\mc D_1(\psi-x)} (x)| \ge |K(\varphi,\theta_{\delta})| - \eps = (1 + \e)^{-d}|K((1+\eps)\varphi,(1+\eps)\theta_{\delta})| - \eps \ge (1 + \e)^{-d}\hspace{-.4cm}\inf_{(\varphi,\psi)\in B} |A(\varphi,\psi)| -\eps.
$$
Since $\eps > 0$ was arbitrary, we conclude the proof.
\enp

\bep[Proof of Lemma \ref{lemma_one_ball_asymptotically_optimal}]
Throughout the proof we rely on the  interpretation of the optimization problem in Lemma \ref{lemma_optimization_problem_for_NNG}.  We set $M\coloneqq c_\ms{max}+1$, and let $(\varphi,\psi)\in B$. Then, we represent $\vp$ as $\vp = \{x_1, \dots, x_m\}$ such that $D_1 \ge D_2 \ge \cdots \ge D_m$, where $D_i \coloneqq \mc D_1(\psi-x_i)$. Next, we  define the normalized $\a$-weighted distances by $\gamma_i \coloneqq {D_i^\a}/(\sum_{j \le m}D_j^\a)$  emphasizing that $D_i\ge \g_i^{1/\a}$ because the denominator is at least 1. For the first part of the lemma, we will distinguish between the two cases that the maximal nearest neighbor distance of a configuration is large or small.
\medskip

\textbf {Case 1: $\boldsymbol{\g_1 \le 1/M}$.}
Note that due to \eqref{FIN} each point in $\R^d$ is contained in at most  $c_\ms{max}$ balls $B_{D_i}(x_i)$, $i \le m$. Thus,
\begin{equation}\label{inequality_union_balls_to_sum}
	|\cup_{i \le m} B_{D_i}(x_i)| \ge \frac1{c_\ms{max}}\sum_{i\le m} |B_{D_i}(x_i)| = \frac{\k_d}{c_\ms{max}} \sum_{i \le m} D_i^d \ge \frac{\k_d}{c_\ms{max}} \sum_{i\le m}\g_i^{d/\a}.
\end{equation}
Now, we formally modify the weights $\{\g_i\}_{i \le m}$ to decrease this sum.  More precisely, we can decrease the values of $\gamma_l$ for $l\in\{M+1,\dots, m\}$ and simultaneously increase some of $\gamma_1,\dots,\gamma_M$ until they are all equal to $1/M$, while keeping $\sum_{i \le m}\g_i =1$. Since concavity implies that $y^{d/\a}+z^{d/\a}\ge (y+z)^{d/\a}$ for $y,z\ge 0$, we deduce that this weight modification only decreases the sum of the $d/\a$-weighted values of the $\gamma_i$'s compared to \eqref{inequality_union_balls_to_sum}. Thus,
\begin{equation}\label{inequality_sum_to_one_ball}
	\frac{\k_d}{c_\ms{max}} \sum_{i \le m} \g_i^{d/\a} \ge \frac{\k_d}{c_\ms{max}} \sum_{i \le M} M^{-d/\a} = \frac{\k_d}{c_\ms{max}} M^{1-d/\a} > \k_d = |B_1(0)|,
\end{equation}
for $\a$ sufficiently large, depending only on $c_\ms{max}$ and $d$.
\medskip

\textbf {Case 2: $\boldsymbol{\g_1  > 1/M}$.} First, we decompose the volume of the union of balls as
$$
|\cup_{i \le m} B_{D_i}(x_i)| = |B_{D_1}(x_1)| + |\cup_{i=2}^{m} \big(B_{D_i}(x_i)\sm B_{D_1}(x_1)\big)|.
$$
Now, note that in the NNG, the balls $B_{D_i}(x_i)$ and $B_{D_1}(x_1)$ cannot fully overlap since $x_i$ cannot be in the interior of $B_{D_1}(x_1)$ and vice versa. Even after subtracting $B_{D_1}(x_1)$, the volume of the remaining shape is still larger than half of its original volume. Thus, by concavity,
 \begin{align*}
|\cup_{i \le m} B_{D_i}(x_i)| - \k_d\g_1^{d/\a}
\ge   \frac1{c_\ms{max}}\sum_{i=2}^{m} \big|B_{D_i}(x_i)\sm B_{D_1}(x_1)\big| \ge  \frac{\k_d}{2c_\ms{max}}\sum_{i=2}^{m} \g_i^{d/\a} \ge  \frac{\k_d}{2c_\ms{max}} (1-\g_1)^{d/\a}.
\end{align*}
Next, since the minimum of a concave function is attained at the boundary,
\begin{equation}\label{inequality_one_ball_optimal_large_ball}
\k_d \g_1^{d/\a} + \frac{\k_d}{2c_\ms{max}} (1-\g_1)^{d/\a} \ge \k_d\min\big\{1,M^{-d/\a} + \f1{2c_\ms{max}} (1-1/M)^{d/\a}\big\} \ge \k_d
\end{equation}
for $\a$ sufficiently large depending on $c_\ms{max}$ and $d$. We summarize the requirements that $\a$ was supposed to be sufficiently large by writing $\a\ge \a_0$ with $\a_0$ depending on $c_\ms{max}$ and $d$. Finally, we point out that the configurations $(\{{0}\},\{{0},(1,0,\dots,0)\})$ are in $B$ since $c_\ms{INF}=2$ for the NNG and it yields the influence zone that is a ball with radius 1 when using the interpretation of the optimization problem for the NNG derived in Lemma \ref{lemma_optimization_problem_for_NNG}. Thus, the volume of the unit ball can indeed be approached by the infimum which gives the first part of Lemma \ref{lemma_one_ball_asymptotically_optimal}.

For the second part, fix $\delta>0$ and let configurations $(\varphi,\psi)\in B$ satisfy $\g_1 < 1-\delta$. We repeat the case distinction that we conducted in the first part and without any adjustments \eqref{inequality_union_balls_to_sum} and \eqref{inequality_sum_to_one_ball} show that if $\g_1 \le 1/M$, there exists an $\eps_1>0$ depending on $c_\ms{max}$ and $d$ such that
$$|\cup_{i \le m} B_{D_i}(x_i)| \ge (1+\eps_1)\k_d$$
for $\a\ge\a_0$.
In the case that $1/M \le \g_1 < 1-\delta$, we can perform a similar calculation as the one that lead to (\ref{inequality_one_ball_optimal_large_ball}) and, by concavity as well as by the fact that the sum of strictly concave functions is again strictly concave, we arrive at
\begin{align*}
|\cup_{i\le m} B_{D_i}(x_i)| &\ge \k_d \g_1^{d/\a} + \frac{\k_d}{2c_\ms{max}} (1-\g_1)^{d/\a} \\
&\ge \k_d\min\big\{ (1-\delta)^{d/\a} + \frac{\delta^{d/\a}}{2c_\ms{max}} ,  M^{-d/\a} + \frac{1}{2c_\ms{max}} (1-1/M)^{d/\a}\big\} \ge (1+\eps_2)\k_d
\end{align*}
for an $\eps_2>0$ depending on $\delta, c_\ms{max}$ and $d$ if $\a\ge\a_0$. Taking $\eps=\min\{\eps_1,\eps_2\}$ concludes the proof. 
\enp
A slightly altered version of the proof of Lemma \ref{lemma_one_ball_asymptotically_optimal} would also work for the undirected NNG. One would have to approximate $(\{0\},\{0,(1,0,\dots,0)\})$ by putting an additional point close to $(1,0,\dots,0)$ to guarantee that the score of the origin is equal to $1$. There are some reasons why the bidirected version does not admit $\kappa_d$ as solution of its optimization problem for large $\a$. First, Lemma \ref{lemma_optimization_problem_for_NNG} does not hold anymore for the bidirected NNG. Another reason is that for $(\{0\},\{0,(1,0,\dots,0)\})$, the value of the score function is $\xi^{(\a)}(\psi)\le 1/2 < 1$ and cannot be approximated with elements of $B$ that yield a score of approximately 1 for the origin while maintaining an influence zone with volume about $\kappa_d$.

\subsection{Theorem \ref{theorem_conditioned_convergence} a) for the (un-/bidirected) kNN and the $\b$-skeleton}

Recall that we need to prove that the optimization problems of the graphs described in Section \ref{sec_application} admit strictly positive solutions. Underlie any of those graphs and let $(\varphi,\psi)\in B$. Note that this implies that
\begin{equation}\label{inequality_positive_sol_optimization}
\sum_{x\in\varphi}\sum_{y\in\mc{E}(\psi-x)} \underbrace{|x-y|^\a}_{\eqqcolon\lambda_{x,y}} \ge \sum_{x\in\varphi} \EE(\psi-x) \ge 1,
\end{equation}
due to the definitions of $\EE$ and $B$ that we recall from \eqref{equation_out_neighbors} and \eqref{definition_sets_optimization}. First, we derive a lower bound for $|A(\varphi,\psi)|$ in terms of a volume of a union of suitable balls. This will be done separately for the (un-/bidirected) kNN and the $\b$-skeletons. After that we can consider both cases simultaneously.

\begin{description}
	\item[(un-/bidirected) kNN:]  First, since $\#\psi\ge c_\ms{INF}=k+1$ we know that for any $x\in\varphi$, an additional node within the interior of $B_{\mc{D}_k(\psi)}(x)$ would delete a vertex in $\{y\co y\in\EE(\psi-x)\}$. The influence zone prohibits such nodes from which we deduce that $|A(\varphi,\psi)| \ge |\cup_xB_{\mc{D}_k(\psi)}(x)| = |\cup_x\cup_y B_{\lambda_{x,y}^{1/\a}} (x)|$. We intentionally let the balls after the equality sign overlap to avoid being forced to distinguish between (un-/bidirected) kNN and $\beta$-skeleton below.
\item[$\bm{\beta}$-skeleton:] For $x\in\varphi$ and $y\in\mc{E}(\psi-x)$, define $h(x,y) \coloneqq (x + y)/2$ as the midpoint between $x$ and $y$. The $\beta$-skeleton for $\beta>1$ is a subgraph of the GG. Therefore, any node put in the ball $B_{\lambda_{x,y}^{1/\a}/2}(h(x,y))$ would removes the edge between $x$ and $y$. Thus, $|A(\varphi,\psi)| \ge |\cup_x\cup_y B_{\lambda_{x,y}^{1/\a}/2} (h(x,y))|.$
\end{description}

Now, enumerate the $\lambda_{x,y}$ decreasingly, i.e., $\lambda_1\ge \lambda_2 \ge \cdots$. Further, we set $\g_i = \lambda_i/(\sum_j \lambda_j)$ achieving that $ \lambda_i\ge \g_i$ due to \eqref{inequality_positive_sol_optimization}. Because of \eqref{FIN} and the bound on the maximal node degree, every point $y\in\R^d$ is contained in at most $(c_\ms{max}+1)^2$ of these balls. Thus,
$$|A(\varphi,\psi)| \ge \sum_i \frac1{(c_\ms{max}+1)^2} |B_{\lambda_i^{1/\a}/2} (0)| \ge \sum_i \frac1{2^d(c_\ms{max}+1)^2} |B_{\g_i^{1/\a}} (0)| = \sum_i \frac{\kappa_d\g_i^{d/\a}}{2^d(c_\ms{max}+1)^2}.$$
Now, as in the proof of Lemma \ref{lemma_one_ball_asymptotically_optimal}, we use concavity to arrive at
$$\sum_i \frac{\kappa_d\g_i^{d/\a}}{2^d(c_\ms{max}+1)^2} \ge \frac{\kappa_d}{2^d(c_\ms{max}+1)^2} \bigg(\sum_i \gamma_i\bigg)^{d/\a} = \frac{\kappa_d}{2^d(c_\ms{max}+1)^2} > 0.$$
Thus, Theorem \ref{theorem_conditioned_convergence} a) becomes applicable.

%
%PART A 
%
\section{Proof of Theorem \ref{theorem_main_upper_tails}}\label{sec_proof_1}

The proof of Theorem \ref{theorem_main_upper_tails} is split up into the upper bound (Section \ref{subsec_up}) and the lower bound (Section \ref{subsec_low}).

\subsection{Upper bound}
\label{subsec_up}
We will follow the strategy that has already been successfully applied in \cite{sphere}, and divide the contributions to $H_n$ into those coming from small or large scores. Then, these are treated separately by the following lemmas, which are shown after the proof of the upper bound of Theorem \ref{theorem_main_upper_tails}. For convenience, we let 
\begin{equation}\label{definition_stabilization_max}
\RR_n(X) \coloneqq \max_{x\in X\cap Q_n}\mc R(X-x)
\end{equation}
denote the maximal stabilization radius in the sampling window, cf.\ \eqref{rr_eq}. We start by bounding summands with small contributions through a Poisson functional concentration inequality from \cite{peccati} to verify that these cannot contribute substantially to the excess.

%
%HOEFFDING
%

\bel[Upper bound for contribution of small summands]
\label{lemma_hoeffding}
Let $\eps \in (0,1)$ and $a\in (0,(1-d/\a)/2)$. Then,
\begin{equation}
\limsup_{n\to\infty} \frac1{n^{d^2/\a}}\log\P\bigg(\frac1{n^d}\sum_{x\in X\cap Q_n} \xi^{(\a)}(X-x)\one\{\xi^{(\a)}(X-x) < n^a\} > \mu_\alpha + \eps r,\RR_n(X) \leq n\bigg) = -\infty.
\end{equation}
\enl

%
%Pois conc
%
Next, we use a concentration result for Binomial random variables from \cite[Lemma 1.1]{poisson_conc} to bound the number 
\begin{equation}
J_n^{(a)}(X)\coloneqq \#\mc J_n^{(a)}(X)\coloneqq\#\{x\in X\cap Q_n \co \xi^{(\a)}(X-x)\ge n^a\}.
\end{equation}
of $x\in X\cap Q_n$ that have a score of at least $n^a$.

\bel[Upper bound for number of large summands]
\label{lemma_poisson_conc}
Let $a\in (0,1)$ and $\eps \in (0,ad/\a)$. Then,
\begin{equation}
\limsup_{n\to\infty} \frac1{n^{d^2/\a}}\log\P\big(J_n^{(a)}(X) > n^{d^2/\a-\eps}\big) = -\infty.
\end{equation}
\enl

%
%UN BOUND
%
Further, we bound the probability that a small number of Poisson points carries a lot of the excess weight.
\bel[Upper bound for condensation probability]
\label{lemma_sum_optimization_problem}
Let $m,n\ge1$ and $\tau>0$. Then, 
\begin{align}
\begin{split}
&\P\bigg(\sum_{x\in\mc J_n^{(a)}(X)} \xi^{(\a)}(X-x) \ge \tau, J_n^{(a)} \le m, \RR_{3n}(X) \le n\bigg) \\
&\le (I_d c_\ms{max} + 1)^4 m^2 (5n)^{d 2 (I_d c_\ms{max} + 1)^2 m} \exp\Big(-\tau^{d/\a} \inf_{(\vp,\psi)\in B} |A(\vp,\psi)|\Big).
\end{split}
\end{align}
\enl
Before proving these lemmas, we apply them to get the upper bound.

%
%PROOFS
%
\bep[Proof of the {upper bound} of Theorem \ref{theorem_main_upper_tails}]
%
%UPPER BOUND
%
Let  $a\in (0,(1-d/\a)/2)$ and $\eps\in(0,ad/\a)$. Then,
\begin{align}\label{inequality_upper_bound_main_ineq}
\begin{split}
&\P\bigg(\sum_{x\in X\cap Q_n} \xi^{(\a)}(X-x) > \mu_\alpha n^d + rn^d\bigg) \\
&\le \P\bigg(\sum_{x\in X\cap Q_n} \xi^{(\a)}(X-x)\one\{\xi^{(\a)}(X-x) < n^a\} - \mu_\alpha n^d > \eps rn^d\bigg) \\
&\quad + \P\bigg(\sum_{x\in X\cap Q_n} \xi^{(\a)}(X-x)\one\{\xi^{(\a)}(X-x) \ge n^a\} \ge (1-\eps)rn^d\bigg) \\
&\le \P\bigg(\frac1{n^d}\sum_{x\in X\cap Q_n} \xi^{(\a)}(X-x)\one\{\xi^{(\a)}(X-x) < n^a\} > \mu_\alpha + \eps r,\RR_n(X) \le n\bigg) + \P\big(J_n^{(a)}(X) > n^{d^2/\a-\eps}\big) \\
&\quad + 2\P(\RR_{3n}(X) > n) + \P\bigg(\sum_{x\in\mc J_n^{(a)} (X)} \xi^{(\a)}(X-x) \ge (1-\eps)rn^d, J_n^{(a)}(X)\le n^{d^2/\a-\eps}, \RR_{3n}(X) \le n\bigg).
\end{split}
\end{align}

From Lemmas \ref{lemma_hoeffding} and \ref{lemma_poisson_conc},  we know that with our choices of $a$ and $\eps$, the first two summands after the last inequality of \eqref{inequality_upper_bound_main_ineq} do not play a role in large volume asymptotics.
Moreover, with the help of Markov's inequality and Mecke's formula \cite[Theorem 4.4]{lastpenrose} we get that
\begin{align}\label{inequality_stabilization_to_cones}
\begin{split}
\P(\RR_{3n}(X) > n) &= \P(\#\{x\in X\cap Q_{3n}\colon \RR(X-x) > n\} \ge 1) \le \E[\#\{x\in X\cap Q_{3n}\colon \RR(X-x) > n\}] \\
&= \E\Big[\sum_{x\in X\cap Q_{3n}} \one\{\RR(X-x) > n\}\Big] = \int_{ Q_{3n}} \P(\RR((X\cup\{x\})-x) \ge n)dx \\
&\le \int_{ Q_{3n}} \sum_{i\le I_d} \P(\SS_i((X\cup\{x\})-x) \ge n)dx.
\end{split}
\end{align}
From here, due to the characteristics of \eqref{STA} it is implied that for each $i\le I_d$ and $r>0$ it holds that
$$|S_i\cap B_r(0)| \ge r^d \underbrace{\min_{j\le I_d} |S_j\cap B_1(0)|}_{ \eqqcolon c_\ms{cones}}$$
and by applying a Poisson concentration bound \cite[Lemma 1.2]{poisson_conc} for a large enough $n$, we can continue our computations for each $i\le I_d$ and $x\in Q_{3n}$ with
\begin{align*}
\P(\SS_i((X\cup\{x\})-x) \ge n) &\le \P\big(X(S_i\cap B_{n/c_\ms{STA}}(0)) \le c_\ms{STA}\big) \\
&\le \exp\bigg(-c_\ms{cones} \Big(\frac{n}{c_\ms{STA}}\Big)^d + c_\ms{STA} - c_\ms{STA}\log\Big(\frac{c_\ms{STA}^{d+1}}{c_\ms{cones} n^d}\Big)\bigg),
\end{align*}
where we recall that if $X$ is interpreted as a Poisson random measure, we can denote the random number of points in a Borel set by $X(\cdot)$. Therefore, continuing from \eqref{inequality_stabilization_to_cones}, we arrive at
\begin{equation}\label{inequality_stabilization_radius}
\frac1{n^{d^2/\a}}\log \P(\RR_{3n}(X) > n) \le - \frac{c_\ms{cones}}{c_\ms{STA}^d} n^{d(1-d/\a)} + \frac{c_\ms{STA}}{n^{d^2/\a}} - \frac1{n^{d^2/\a}}\log\Big(\frac{c_\ms{STA}^{d+1}}{c_\ms{cones} n^d}\Big) \overset{n\to\infty}{\longrightarrow} -\infty.
\end{equation}
Thus, it remains to consider the fourth summand after the last inequality of \eqref{inequality_upper_bound_main_ineq}. Here, Lemma \ref{lemma_sum_optimization_problem} yields
\begin{align*}
	&\limsup_{n \to \ff}\frac1{n^{d^2/\a}} \log\P\bigg(\sum_{x\in\mc J_n^{(a)} (X)} \xi^{(\a)}(X-x) \ge (1-\eps)rn^d, J_n^{(a)}(X)\le n^{d^2/\a-\eps}, \RR_{3n}(X) \le n\bigg) \\
&\le -((1-\eps)r)^{d/\a} \inf_{(\vp,\psi)\in B} |A(\vp,\psi)|.
\end{align*}
In brief, we arrive at
$$\limsup_{n \to \ff}\frac1{n^{d^2/\a}}\log\P(H_n > \mu_\a + r) \le -((1-\eps)r)^{d/\a} \inf_{(\vp,\psi)\in B} |A(\vp,\psi)|.$$
Letting $\eps\downarrow 0$ concludes the proof of the upper bound.
\enp

In the rest of this subsection, we will prove Lemmas \ref{lemma_hoeffding}, \ref{lemma_poisson_conc} and \ref{lemma_sum_optimization_problem}. The essential ingredient for the proof of Lemma \ref{lemma_hoeffding} is a concentration bound from \cite[Corollary 3.3 i)]{peccati}.
\bep[Proof of Lemma \ref{lemma_hoeffding}]
We start by introducing some of the notation from \cite{peccati}. For the Poisson process $X$, we define the functional
\begin{equation}\label{definition_functional_conc}
F_n^{(\a)}(X) \coloneqq \sum_{x\in X\cap Q_n} \xi^{(\a)}(X\cap Q_{3n}-x)\one\{\xi^{(\a)}(X\cap Q_{3n}-x) < n^a\}.
\end{equation}
Before we can apply the concentration bound, we need to find a link between the typical value $n^d \mu_\a$ and the expectation of the functional defined in \eqref{definition_functional_conc}. We can find the connection using that, because of \eqref{STA}, under the event $\{\RR_n(X) \leq n\}$, this functional is equal to the one considered in Lemma \ref{lemma_hoeffding}
\begin{align*}
n^d\mu_\a &\ge \E\bigg[\sum_{x\in X\cap Q_n}\xi^{(\a)}(X-x)\one_{\{\xi^{(\a)}(X-x) < n^a\}}\bigg] \ge \E\bigg[\sum_{x\in X\cap Q_n}\xi^{(\a)}(X-x)\one_{\{\xi^{(\a)}(X-x) < n^a\}}\one_{\{\RR_n(X) \le n\}}\bigg] \\
&= \E\bigg[\sum_{x\in X\cap Q_n}\xi^{(\a)}(X\cap Q_{3n}-x)\one{\{\xi^{(\a)}(X\cap Q_{3n}-x) < n^a, \, \RR_n(X) \le n\}}\bigg] \\
&= \E\bigg[\sum_{x\in X\cap Q_n}\xi^{(\a)}(X\cap Q_{3n}-x)\one{\{\xi^{(\a)}(X\cap Q_{3n}-x) < n^a\}}(1-\one{\{\RR_n(X) > n\}})\bigg] \\
&= \E[F_n^{(\a)}(X)] - \E\bigg[\sum_{x\in X\cap Q_n}\xi^{(\a)}(X\cap Q_{3n}-x)\one{\{\xi^{(\a)}(X\cap Q_{3n}-x) < n^a,\,\RR_n(X) > n\}}\bigg] \\
&\ge \E[F_n^{(\a)}(X)] - \E[X(Q_n)n^a\one{\{\RR_n(X) > n\}}].
\end{align*}
Subsequently, the Cauchy-Schwarz inequality yields
$$n^d\mu_\a \ge \E[F_n^{(\a)}(X)] - \sqrt{\E[X(Q_n)^2 n^{2a}] \E[\one_{\{\RR_n(X) > n\}}]} =  \E[F_n^{(\a)}(X)] - \sqrt{(n^{2d}+n^d)n^{2a} \P(\RR_n(X)>n)}.$$
As argued in \eqref{inequality_stabilization_radius}, the factor $\P(\RR_n(X)>n)$ decays exponentially with $n$ and consequently, we can assume that $n$ is large enough to guarantee that  $n^d\mu_\a \ge \E[F_n^{(\a)}(X)] - \eps$.  Therefore,
\begin{align}\label{inequality_functional_concentration_application}
\begin{split}
&\P\bigg(\sum_{x\in X\cap Q_n} \xi^{(\a)}(X-x)\one_{\{\xi^{(\a)}(X-x) < n^a\}} > n^d\mu_\alpha + n^d\eps r,\RR_n(X) \leq n\bigg) \le \P(F_n^{(\a)}(X) > n^d\mu_\a +n^d\eps r) \\
&\le \P(F_n^{(\a)}(X) > \E[F_n^{(\a)}(X)] + n^d\eps r - \eps).
\end{split}
\end{align}
Further, we need the difference operator $D_y$, $y \in \R^d$, defined by
$D_yF_n^{(\a)}(X) \coloneqq F_n^{(\a)}(X\cup\{y\}) - F_n^{(\a)}(X).$
For $\beta\ge0$, we now set
\begin{equation}\label{definition_V_conc}
V_\beta^+(F_n^{(\a)}(X)) \coloneqq \int_{\R^d} \big(D_yF_n^{(\a)}(X)\one_{\{D_yF_n^{(\a)}(X) \le \beta\}}\big)^2 dy + \sum_{x\in X} \big(D_xF_n^{(\a)}(X\sm\{x\})\one_{\{D_xF_n^{(\a)}(X\sm\{x\}) > \beta\}}\big)^2.
\end{equation}
To apply \cite[Corollary 3.3 i)]{peccati}, we need to find an almost sure upper bound for $V_\beta^+(F_n^{(\a)})$. Points outside of $Q_{3n}$ do not affect the functional. Thus, choosing $y\in\R^d\sm Q_{3n}$ in the difference operator has no effect and yields $D_yF_n^{(\a)}(X) = 0$. Besides, due to \eqref{FIN}, adding a point to any configuration can only affect the outgoing edges of $c_\ms{max}$ nodes and the degree of each node is bounded by $c_\ms{max}$ as well. Hence, $\sup_{y \in Q_{3n}}|D_yF_n^{(\a)}(X)| \le(c_\ms{max}+1)^2 n^a \eqqcolon \beta$, and by the same reasoning, $\sup_{x \in X}|D_xF_n^{(\a)}(X\sm\{x\})| \le \beta$. Thus, we bound \eqref{definition_V_conc} by
$$V_\beta^+(F_n^{(\a)}(X)) \le \int_{Q_{3n}} (c_\ms{max}+1)^4 n^{2a} dy = ((c_\ms{max}+1)^2 n^a)^2(3n)^d.$$
Then, by applying \cite[Corollary 3.3 i)]{peccati},
\begin{align*}
\P(F_n^{(\a)}(X) > \E [F_n^{(\a)}(X)] + n^d\eps r - \eps) &\le \exp\Big(-\frac{n^d\eps r - \eps}{2|\beta|}\log\Big(1+\frac{|\beta|(n^d\eps r - \eps)}{((c_\ms{max}+1)^2 n^a)^2(3n)^d}\Big)\Big) \\
&= \exp\bigg(-\frac{n^d\eps r - \eps}{2(c_\ms{max}+1)^2 n^a}\log\Big(1+\frac{(c_\ms{max}+1)^2 n^a (n^d\eps r - \eps)}{(c_\ms{max}+1)^4 n^{2a}(3n)^d}\Big)\bigg) \\
&= \exp\bigg(-\frac{n^{d-a}(\eps r - \eps/n^d)}{2(c_\ms{max}+1)^2}\log\Big(1+\frac{n^{-a}(\eps r - \eps/n^d)}{(c_\ms{max}+1)^2 3^d}\Big)\bigg)
\end{align*}
if $n^d\eps r - \eps \ge 0$. Finally, with the help of \eqref{inequality_functional_concentration_application}
$$\limsup_{n\to\infty} \frac1{n^{d^2/\a}} \log \P\bigg(\sum_{x\in X\cap Q_n} \xi^{(\a)}(X-x)\one_{\{\xi^{(\a)}(X-x) < n^a\}} > n^d\mu_\a +n^d\eps r, \RR_n(X) \le n\bigg) = - \infty$$
for $a\in(0,(1-d/\a)/2)$.
\enp

\begin{figure}[H]
\centering
\begin{tikzpicture}[scale=0.8]
%Big Cube [0,3]^3
\draw (3,3,0)--(0,3,0)--(0,3,3)--(3,3,3)--(3,3,0)--(3,0,0)--(3,0,3)--(0,0,3)--(0,3,3);
\draw (3,3,3)--(3,0,3);
\draw[gray!50] (3,0,0)--(0,0,0)--(0,3,0);
\draw[gray!50] (0,0,0)--(0,0,3);

\foreach \x in{1,...,3}
{
\foreach \y in{1,...,3}
{
\foreach \z in{3}
{   
\pgfmathsetmacro\k{(\x+\y*\y/9+\z*\z*\z/27)*10};
\draw[fill=blue!\k] ({\x*0.25},{(\y-1)*0.25},{\z*0.25})--({(\x-1)*0.25},{(\y-1)*0.25},{\z*0.25})--({(\x-1)*0.25},{\y*0.25},{\z*0.25})--({\x*0.25},{\y*0.25},{\z*0.25})--({\x*0.25},{(\y-1)*0.25},{\z*0.25});

}
}
}

\foreach \x in{3}
{
\foreach \y in{1,...,3}
{
\foreach \z in{1,...,3}
{   
\pgfmathsetmacro\k{(\x+\y*\y/9+\z*\z*\z/27)*10};
\draw[fill=blue!\k] ({\x*0.25},{(\y-1)*0.25},{(\z-1)*0.25})--({\x*0.25},{(\y-1)*0.25},{\z*0.25})--({\x*0.25},{\y*0.25},{\z*0.25})--({\x*0.25},{\y*0.25},{(\z-1)*0.25})--({\x*0.25},{(\y-1)*0.25},{(\z-1)*0.25});

}
}
}

\foreach \x in{1,...,3}
{
\foreach \y in{3}
{
\foreach \z in{1,...,3}
{
\pgfmathsetmacro\k{(\x+\y*\y/9+\z*\z*\z/27)*10};
\draw[fill=blue!\k] ({(\x-1)*0.25},{\y*0.25},{(\z-1)*0.25})--({(\x-1)*0.25},{\y*0.25},{\z*0.25})--({\x*0.25},{\y*0.25},{\z*0.25})--({\x*0.25},{\y*0.25},{(\z-1)*0.25})--({(\x-1)*0.25},{\y*0.25},{(\z-1)*0.25});

}
}
}

\end{tikzpicture}
\hspace{3cm}
\begin{tikzpicture}[scale=0.9]
\draw (0,0) rectangle (3.7,3.7);
\draw (0,0) rectangle (0.25,0.25);
\draw[pattern=north east lines] (0,0.25) rectangle (0.25,0.5);
\draw (0,0.75) rectangle (0.25,1);
\draw (0,0.5)[fill=black] rectangle (0.25,0.75);
\draw[fill=gray!50] (0.25,0) rectangle (0.5,0.25);
\draw (0.25,0.25) rectangle (0.5,0.5);
\draw (0.25,0.5) -- (0.5,0.25);
\draw[fill=gray!50] (0.25,0.75) rectangle (0.5,1);
\draw[pattern=north west lines] (0.25,0.5) rectangle (0.5,0.75);
\draw (0.75,0) rectangle (1,0.25);
\draw (0.5,0)[pattern=crosshatch] rectangle (0.75,0.25);
\draw[pattern=north east lines] (0.75,0.25) rectangle (1,0.5);
\draw[fill=black] (0.5,0.25) rectangle (0.75,0.5);
\draw[white] (0.5,0.25) -- (0.75,0.5);
\draw (0.5,0.25) rectangle (0.75,0.5);
\draw (0.75,0.75) rectangle (1,1);
\draw (0.75,0.5)[fill=black] rectangle (1,0.75);
\draw (0.5,0.75)[pattern=crosshatch] rectangle (0.75,1);
\draw (0.5,0.5)[fill=gray!70!black] rectangle (0.75,0.75);
\end{tikzpicture}
\caption{Labeling of the boxes in 3D where 27 labels are sufficient and 2D where 9 are sufficient}\label{fig:cubes}
\end{figure}
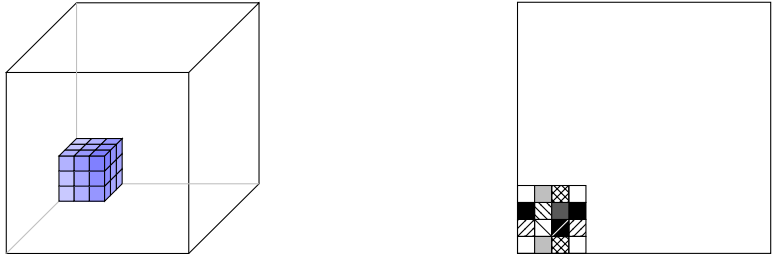

\bep[Proof of Lemma \ref{lemma_poisson_conc}]
Let $a\in(0,1)$. We divide $Q_n$ into a grid consisting of $\lfloor n^{1-a/\a}\rfloor^d$ smaller boxes with side length $l_n \coloneqq n/\lfloor n^{1-a/\a}\rfloor$. The set of all of these cubes is
$$\mc Q \coloneqq \big\{Q \colon Q = l_n z + [-n/2,-n/2+l_n]^d, z\in\{0,\dots,\lfloor n^{1-a/\a}\rfloor - 1\}^d\big\}.$$
Furthermore, we label each box in a way that between two boxes of the same label there are always two boxes with a different label.  For instance, we can label the boxes according to elements of the set $\LL = \{0,1,2\}^d$, thus using $\#\LL = 3^d$ different labels, see Figure \ref{fig:cubes}. For $m\in \LL$, we denote the set of label $m$ cubes by
$$\mc Q^{(m)} \coloneqq \big\{Q \colon Q = l_n z + [-n/2,-n/2+l_n]^d, z=(z_1,\dots,z_d)\in\{0,\dots,\lfloor n^{1-a/\a}\rfloor - 1\}^d \text{ with } z_i\bmod 3 = m_i\big\},$$
so that
$\#\mc Q^{(m)} = \lfloor n^{1-a/\a}\rfloor^d/3^d + o( n^{1-a/\a}) \eqqcolon K_n.$

Setting $P_n\coloneqq (n^a/c_\ms{max})^{1/\a}$, we start bounding the considered probabilities using the bounded node degree:
\begin{align}\label{inequality_subcubes1}
\begin{split}
\P(J_n^{(a)}(X) > n^{d^2/\a-\eps}) &\le \P\big(\xi^{(\a)}(X-x) > n^a \text{ for all }x\text{ in some }\vp\subseteq X\cap Q_n\text{ with }\#\vp\ge n^{d^2/\a-\eps}\big) \\
&\le \P\big(\max_{y\in\mc E(X-x)} |y| > P_n \text{ for all }x\text{ in some }\vp\subseteq X\cap Q_n\text{ with }\#\vp\ge  n^{d^2/\a-\eps}\big).
\end{split}
\end{align}
We now thin out the configuration consisting of all $x$ as in the previous line as follows. Starting with any point $x\in\vp$, we omit all points of $\vp$ that are at distance at most $P_n$ to $x$. According to \eqref{FIN2} with $M=P_n$, this operation removes at most $c_\ms{max}-1$ points. Repeating iteratively for the other points of $\vp$ yields a configuration $\vp$ that contains at least $N_n \coloneqq n^{d^2/\a-\eps}/c_\ms{max}$ nodes satisfying $\max_{y\in\mc E(X-x)} |y| > P_n$ and $|x-y| > P_n$ for all $x,y\in\vp$ with $x\ne y$. Thus, we can continue in \eqref{inequality_subcubes1} with
\begin{align}\label{inequality_subcubes2}
\begin{split}
&\P\big(\max_{y\in\mc E(X-x)} |y| > P_n \text{ for all }x\text{ in some }\vp\subseteq X\cap Q_n\text{ with }\#\vp\ge n^{d^2/\a-\eps}\big) \\
&\le \P\big(\max_{y\in\mc E(X-x)} |y| > P_n \text{ and } |x-y| > P_n \text{ for all }x\ne y\text{ in some }\vp\subseteq X\cap Q_n\text{ with }\#\vp\ge N_n \big).
\end{split}
\end{align}
Next, we note that in a ball of radius $\sqrt dl_n$, only a limited number of points can be placed such that all of their mutual distances are larger than $P_n$. For large $n$, this number is bounded by the number of balls with radius $(n^a/c_\ms{max})^{1/\a}/2$ that fit in a ball with radius $4\sqrt d n^{a/\a}$ such that none of the smaller balls overlap. The fraction of the volume of $B_{4\sqrt d n^{a/\a}}(0)$ and the volume of $B_{(n^a/c_\ms{max})^{1/\a}/2}(0)$ yields the bound $8^d c_\ms{max}^{d/\a} d^{d/2}$ for $n$ large. Thus, after setting $M_n \coloneqq N_n/(8^d c_\ms{max}^{d/\a + 1} d^{d/2})$, we can use this argument to proceed in \eqref{inequality_subcubes2} and estimate for $n$ sufficiently large
\begin{align}\label{inequality_subcubes3}
\begin{split}
&\P\big(\max_{y\in\mc E(X-x)} |y| > P_n \text{ and } |x-y| > P_n \text{ for all }x\ne y\text{ in some }\vp\subseteq X\cap Q_n\text{ with }\#\vp\ge N_n \big) \\
&\le \P\big(\max_{y\in\mc E(X-x)} |y| > P_n \text{ and } |x-y| > \sqrt{d} l_n \text{ for all }x\ne y\text{ in some }\vp\subseteq X\cap Q_n\text{ with }\#\vp\ge M_n \big).
\end{split}
\end{align}
In the event on the right-hand side of \eqref{inequality_subcubes3}, each hypercube $Q\in\mc Q$ contains at most one node that has an edge larger than $P_n$. Further if $\max_{y\in\mc E(X-x)} |y| > P_n$ holds for an $x\in X$,  then \eqref{STA} gives that $\RR(X-x) > P_n$. Thus, by a union bound, we arrive at
\begin{align}\label{inequality_subcubes4}
\begin{split}
&\P\big(\max_{y\in\mc E(X-x)} |y| > P_n \text{ and } |x-y| > \sqrt{d} l_n \text{ for all }x\ne y\text{ in some }\vp\subseteq X\cap Q_n\text{ with }\#\vp\ge M_n \big) \\
&\le \sum_{m\in \LL} \P\big(\#\{Q\in \mc Q^{(m)} \co \max_{y\in\mc E(X-x)} |y| > P_n \text{ for some } x\in Q\cap X\} \ge M_n/3^d\big) \\
&\le \sum_{m\in \LL} \P\big(\#\{Q\in \mc Q^{(m)} \colon \RR(X-x) \ge P_n \text{ for some } x\in Q\cap X\} \ge M_n/3^d\big).
\end{split}
\end{align}
With a calculation performed in the same fashion as in \eqref{inequality_stabilization_radius} we get
$$\P\big(\max_{x \in Q \cap X}\RR(X-x) \ge P_n \big) \le \int_Q\P\big(\RR(X\cup\{x\}-x) \ge (n^a/c_\ms{max})^{1/\a}\big)dx \le l_n^ d e^{-cn^{ad/\a}}$$
for $n$ large enough and a value $c>0$. Next, note that $l_n\ge P_n$. Thus, for a fixed $m\in\LL$ the events of finding a Poisson point with a stabilization radius exceeding $P_n$ in a box $Q$ are independent for different choices of $Q \in \mc Q^{(m)}$. Therefore, a binomial concentration bound \cite[Lemma 1.1]{poisson_conc} gives that for each $m\in \LL$
\begin{align}\label{inequality_subcubes5}
\begin{split}
	&\P\big(\#\{Q\in \mc Q^{(m)} \colon \max_{x \in Q \cap X}\RR(X-x) \ge P_n \} \ge M_n/3^d\big) \\
&\le \exp\bigg(-\frac{M_n}{3^d 2} \log\bigg(\frac{M_n/3^d}{K_n l_n^d e^{-cn^{ad/\a}}}\bigg)\bigg)
\end{split}
\end{align}
assuming $n$ is sufficiently large. Now, note that
\begin{equation}\label{inequality_subcubes6}
\lim_{n\to\infty} -\frac1{n^{d^2/\a}}\frac{M_n}{3^d 2}\log\bigg(\frac{M_n/3^d}{K_n l_n^ d e^{-cn^{ad/\a}}}\bigg) = -\infty
\end{equation}
holds if $\eps\in (0,ad/\a)$. Finally, combining \eqref{inequality_subcubes1}, \eqref{inequality_subcubes1}, \eqref{inequality_subcubes2}, \eqref{inequality_subcubes3}, \eqref{inequality_subcubes4}, \eqref{inequality_subcubes5} and \eqref{inequality_subcubes6} yields the desired result.
\enp

\bep[Proof of Lemma \ref{lemma_sum_optimization_problem}]
Let $m\ge1$ and $\tau>0$ and let us assume that we are under the event that we would like to bound in Lemma \ref{lemma_sum_optimization_problem}, Note that due to \eqref{equation_weaker_stabilization_req}, under $\{\RR_{3n}(X) \le n\}$ we have that $\mc E(X-x) = \mc E(X\cap Q_{5n}-x)$ for all $x\in X\cap Q_{3n}$. Under the event $\{J_n^{(a)}(X) \le m\}$, we choose $\vp=\mc J_n^a(X)$ and $\psi'=\cup_{x\in\vp}\mc E_x(X)$. From \eqref{STA}, we obtain configurations $\theta_x$, $x\in\vp\cup\psi'$ with
$\mc E(X-x)=\mc E(\theta_x-x)$
and $\#\theta_x\le I_d c_\ms{STA}$. \eqref{STA} also implies that $\mc E(X-x)=\mc E(\psi-x)$   for every $x\in\vp\cup\psi'$ where $\psi \coloneqq \vp\cup\psi'\cup_{x\in\vp\cup\psi'} \theta_x$.
Note that due to the bound on the stabilization radius, the set $\psi$ is entirely contained in $Q_{5n}$. Moreover, the bounded node degree implies that
\begin{equation}\label{inequality_size_configurations}
\#\psi \le \#\vp + \#\psi' + (\#\vp + \#\psi') I_d c_\ms{STA} \le (I_d c_\ms{STA} + 1) (c_\ms{max} + 1) m \le (I_d c_\ms{max} + 1)^2 m.
\end{equation}
Below, in the case that $\#\psi<c_\ms{INF}$, we add $c_\ms{INF}-\#\psi$ points in $X\cap Q_{5n}$ to $\psi$ to be able to apply \eqref{INF}. To justify that $X(Q_{5n}) \ge c_\ms{INF}-\#\psi$ can be assumed here, we remark that under $\big\{\sum_{x\in\mc J_n^{(a)}(X)} \xi^{(\a)}(X-x) \ge \tau, \RR_{3n}(X) \le n\big\}$ it has to hold that $X(Q_{5n})\ge c_\ms{INF}$. The reason for this is that $\xi^{(\a)}(X-x) > 0$ for some $x\in X\cap Q_n$ implies that $X\cap Q_n$ cannot be empty and $\RR_{3n}(X) \le n$ implies that there have be at least $I_d (c_\ms{max}-1)$ other Poisson points within distance $n$ of any point in $Q_n$. Thus, we even get that $X(Q_{5n})\ge I_d(c_\ms{max}-1) + 1 \ge c_\ms{INF}$, which concludes this argument. Next, using \eqref{inequality_size_configurations} together with $\mc E(X-x)=\mc E(\psi-x)$ and \eqref{INF} we obtain that
\begin{align*}
&\P\bigg(\sum_{x\in\mc J_n^{(a)}(X)} \xi^{(\a)}(X-x) \ge \tau, J_n^{(a)}(X) \le m, \RR_{3n}(X) \le n\bigg) \\
&\le\P\bigg(\sum_{x\in\vp} \xi^{(\a)}(\psi-x) \ge \tau, \text{ for some }\vp\subseteq X\cap Q_n, \#\vp \le m \text{ and } \vp\subseteq\psi\subseteq X\cap Q_{5n}, \\
&\qquad\quad c_\ms{INF} \le \#\psi\le (I_d c_\ms{max} + 1)^2 m, \mc E(\psi-x)=\mc E(X-x) \text{ for all } x\in\vp \cup \cup_{z\in\vp} \mc E_z(\psi)\bigg) \\
&\le\P\bigg(\sum_{x\in\vp} \xi^{(\a)}(\psi-x) \ge \tau, \text{ for some }\vp\subseteq X\cap Q_n, \#\vp \le m \text{ and } \vp\subseteq\psi\subseteq X\cap Q_{5n}, \\
&\qquad\quad c_\ms{INF} \le \#\psi\le (I_d c_\ms{max} + 1)^2 m, \mc E(\psi-x)\subseteq\mc E((\psi\cup\{y\}-x) \text{ for all } y\in X \text{ and }x\in\vp \cup \cup_{z\in\vp}\mc E_z(\psi)\bigg) \\
&=\P\bigg(\sum_{x\in\vp} \xi^{(\a)}(\psi-x) \ge \tau, \text{ for some }\vp\subseteq X\cap Q_n, \#\vp \le m \text{ and } \vp\subseteq\psi\subseteq X\cap Q_{5n}, \\
&\qquad\quad c_\ms{INF} \le \#\psi\le (I_d c_\ms{max} + 1)^2 m, X\cap A(\vp,\psi) = \emptyset\bigg) \eqqcolon (\star).
\end{align*}
We remind the reader of the sets $D_l'$, $l\in\N$ and $\NN'$ that have been defined in Section \ref{sec_model} before Equation \eqref{definition_sets_optimization}. Note that due to the assumptions in \eqref{CON},
$$0=|N_{l+1}| = \int_{\{(x_1,\dots,x_l)\in\R^{dl}:\text{ pw.~distinct}\}} |D(\{x_1,\dots,x_l\})| d(x_1,\dots,x_l)$$
which implies that $|{D}_l'|=0$ and thus, $\NN'$ is a zeroset. In the following let $\mathbf{x}=(x_1,\dots,x_{l_1})$ and $\mathbf{y}=(y_1,\dots,y_{l_2})$ represent $\vp$ and $\psi\sm\vp$, respectively. We will abuse notation and allow $\mathbf{x}$ and $\mathbf{y}$ to be treated as sets. A combination of the union bound, Markov's inequality and Mecke's formula yields
\begin{align*}
(\star) &\le \sum_{0\le l_1,l_2\le (I_d c_\ms{max} + 1)^2 m} \int_{Q_{5n}^{l_2}}\int_{Q_{5n}^{l_1}} \P\bigg(\sum_{x\in\mathbf{x}} \xi^{(\a)}(\mathbf{x}\cup\mathbf{y}-x) \ge \tau,\#\mathbf{y} \ge c_\ms{INF}, X\cap A(\mathbf{x},\mathbf{x}\cup\mathbf{y})=\emptyset\bigg) d\mathbf{x} d\mathbf{y} \\
&= \sum_{0\le l_1,l_2\le (I_d c_\ms{max} + 1)^2 m} \int_{Q_{5n}^{l_2}}\int_{Q_{5n}^{l_1}} \one_{\{\sum_{x\in\mathbf{x}} \xi^{(\a)}(\mathbf{x}\cup\mathbf{y}-x) \ge \tau\}} \one_{\{\mathbf{x}\cup\mathbf{y}\not\in\mc{N}',\,\#\mathbf{y} \ge c_\ms{INF}\}} \exp(-|A(\mathbf{x},\mathbf{x}\cup\mathbf{y})|) d\mathbf{x} d\mathbf{y} \\
&\le (I_d c_\ms{max} + 1)^4 m^2 (5n)^{d 2 (I_d c_\ms{max} + 1)^2 m}\exp\Big(-\tau^{d/\a}\inf_{(\mathbf{x},\mathbf{x}\cup\mathbf{y})\in B} |A(\mathbf{x},\mathbf{x}\cup\mathbf{y})|\Big)
\end{align*}
from which the assertion follows.
\enp

%
%LOWER BOUND
%
\subsection{Lower bound}
\label{subsec_low}

\begin{figure}[H]
\centering
\begin{tikzpicture}[scale=.7]
\draw (0,0) rectangle (10,10);

\draw[pattern=north east lines, pattern color=blue!40] (0,0) rectangle (6.75,6.75);

\input{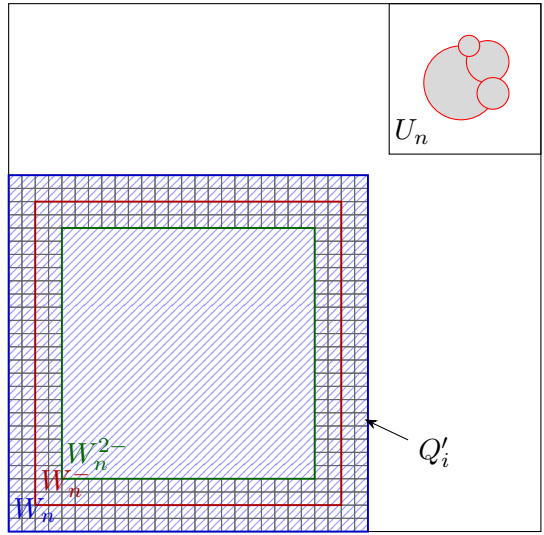}

\draw[line width=0.25mm, green!40!black] (1,1) rectangle (5.75,5.75) node[pos=.14, color=green!40!black, yshift = -0.14cm] {$W_n^{2-}$};
\draw[line width=0.25mm, red!70!black] (0.5,0.5) rectangle (6.25,6.25) node[pos=.1, color=red!70!black, yshift = -0.1cm] {$W_n^-$};
\draw[line width=0.25mm, blue!80!black] (0,0) rectangle (6.75,6.75) node[pos=.06, color=blue!80!black,  xshift = 0.06cm] {$W_n$};
\draw (7.15,7.15) rectangle (10,10) node[pos=.15] {$U_n$};

\node (square_label1) at (8, 1.5) {$Q'_i$};
\draw[-{Stealth}] (square_label1) -- (6.6875, 2.125);

%\node (square_label2) at (7.75, 4.5) {$Q''_i$};
%\draw[-{Stealth}] (square_label2) -- (6.1875, 5.125);

\draw[color=red, fill=gray!30] (8.5,8.5) circle (.7cm);
\draw[color=red, fill=gray!30] (9,8.9) circle (.4cm);
\draw[color=red, fill=gray!30] (9.1,8.3) circle (.3cm);
\draw[color=red, fill=gray!30] (8.65,9.2) circle (.2cm);
\end{tikzpicture}
	\caption{Sketch of $U_n$, $W_n$, $W_n^-$ and $W_n^{2-}$.}\label{fig:sets}
\end{figure}

First, if $\inf_{(\varphi,\psi)\in B} |A(\varphi,\psi)|=\ff$ there is nothing to prove. Thus, throughout the proof of the lower bound we assume that $\inf_{(\varphi,\psi)\in B} |A(\varphi,\psi)|<\ff$. Recall that $\lceil\cdot\rceil$ denotes the ceiling function given by $\lceil t\rceil \coloneqq \min\{m\in\Z\colon m\ge t\}$ for $t\in\R$. The rough idea for the proof of the lower bound is to use separated boxes
$$W_n \coloneqq \big[0,\underbrace{\lceil n-n^{d/\a}\log n- (\log n)^2\rceil}_{\eqqcolon b_n}\big]^d$$
and
$U_n \coloneqq [n-n^{d/\a}\log n, n]^d$
and place the configuration responsible for the excess weight entirely in $U_n$ while letting $W_n$ be responsible for the typical value. The separation is achieved by conditioning on points being close to the boundary of $W_n$. In particular, we introduce a smaller box
$$W_n^{2-} \coloneqq \big[2(\log n)^2, b_n - 2(\log n)^2\big]$$
inside of $W_n$ and condition on a certain amount of points laying in $W_n\sm W_n^{2-}$. This is realized by covering that volume with layers of boxes with side lengths between $\log n$ and $2\log n$, preferably hypercubes with length $\log n$ as pointed out in Figure \ref{fig:sets}. Hence, each box has a volume between $(\log n)^d$ and $(2\log n)^d$. 
%\begin{wrapfigure}{r}{0.42549\textwidth}
%\end{wrapfigure}

Hence, for sufficiently large  $n$ we need at most
\begin{equation}\label{number_of_boxes_to_cover_cube}
\left\lceil \frac{b_n^d-(b_n - 4(\log n)^2)^d}{(\log n)^d}\right\rceil \le \left\lceil\frac{n^{d-1}}{(\log n)^{d-3}}\right\rceil
\end{equation}
additional boxes to cover the space $W_n\sm W_n^{2-}$ entirely. We denote these boxes by $(Q'_i)_i$ and define the event
$$E_n^\ms{good} \coloneqq \big\{X(Q'_i) \in [c_\ms{max},(\log n)^{2d}) \text{ for all } i\big\}$$
that will generate an independence between the functional of Poisson points in $W_n$ and Poisson points in $U_n$. Besides that, we introduce the abbreviation
$$H_n (A,B) \coloneqq \frac1{n^d}\sum_{x\in X\cap A} \xi^{(\a)}(X\cap B-x)$$
for $A,B\subseteq\R^d$. For $\eps < \mu_\a$, we also define the event
\begin{equation}\label{definition_G1}
G_{1,n} \coloneqq \{H_n(W_n,W_n) > \mu_\a - \eps/2\} \cap E_n^\ms{good}.
\end{equation}
The next lemma gives a lower bound for the probability of this event.

	\bel[Lower bound for $\P(G_{1, n})$]\label{lemma_weak_law_of_large_numbers}
It holds that
$\liminf_{n\to\infty} {n^{-d^2/\a}} \log \P(G_{1,n}) \ge 0.$
\enl

We now focus our attention on what happens within $U_n$. We will rescale a configuration in a way such that it is responsible for the entire excess weight and such that there is also enough flexibility to embed the points in open balls to get a configuration that can be attained with positive probability. For the chosen $\eps$, we will use
$$\tau_n \coloneqq ((r+\eps) (1+\eps) n^d)^{1/\a}$$
as parameter for the rescaling. The following lemma will be used to find the proper configuration within $U_n$ to rescale.
\bel[Approximation of optimal configurations]\label{lemma_approximation_optimal_configurations}
	Let $\e > 0$ and  $(\vp, \psi) \in B$. Then, there exists $\de \in (0,1)$ such that the following inequalities hold
\begin{enumerate}
\item[a)] 
$$
|\underbrace{\bigcup_{(z_y)_{y\in\psi}\subseteq B_\delta(0)} A(\{x+z_x \co x\in\vp\},\{y+z_y \co y\in\psi\})}_{\eqqcolon A_{\delta}(\vp,\psi)}|
		\le |A(\vp, \psi)| + \e 
$$
and
\item[b)]
$$\inf_{(z_x)_{x\in\psi} \subseteq B_\delta(0)} \sum_{x\in\vp} \xi^{(\a)}(\{y+z_y \colon y\in\psi\}-(x+z_x)) > 1/(1+\eps).$$
\end{enumerate}
\enl

We insert another lemma to deal with the diameter of the influence zone
\bel[Diameter of bounded influence zone]\label{lemma_diam_influence_zone}
Let $(\varphi,\psi)\in B$ with $|A(\varphi,\psi)| < \ff$. Then, there is $\de \in (0, 1)$ such that $\diameter(A_\de(\varphi,\psi))<\ff$.
\enl

Note that, if we pick $(\varphi,\psi)\in B$ such that $|A(\varphi,\psi)|<\infty$, then, for $\de $ small enough, by Lemma \ref{lemma_diam_influence_zone} the diameter of $\t_nA_{\delta}(\vp,\psi)$ is of order $n^{d/\a}$, while $U_n$ has side length $n^{d/\a}\log n$. This means, we can choose $n$ large enough for $U_n$ to contain a shifted copy of $\t_n A_{\delta}(\vp,\psi)$. Thus, from now on we can assume that $\t_nA_{\delta}(\vp,\psi)$ and also $\cup_{x\in\t_n\psi} B_1(x)$, is entirely contained in $U_n$ if $|A(\varphi,\psi)|<\infty$.

We set
$$A_{\delta, n}^- \coloneqq \t_nA_{\delta}(\vp,\psi) \sm \bigcup_{x\in\tau_n\psi} B_1(x)$$ and similarly to \eqref{definition_G1}, we define the event
\begin{equation}\label{definition_G2}
G_{2,n}(\delta) \coloneqq \big\{X(B_1(x)) = 1 \text{ for all } x\in\tau_n\psi, X(A_{\delta, n}^-) = 0\big\}.
\end{equation}
A bound for its probability is given in the following lemma.

\bel[Lower bound for $\P(G_{2, n}(\de))$]\label{lemma_probability_points_responsible_for_excess_weight}
Let $\de \in (0, 1)$ and $(\vp, \psi) \in B$. Then, 
$$\liminf_{n\to\infty} \frac1{n^{d^2/\a}} \log \P(G_{2,n}(\delta)) \ge -\big( |A(\vp,\psi)|+\eps\big)(r+\eps)^{d/\a}(1+\eps)^{d/\a}.$$
\enl

Now, we can state the proof of the lower bound.

\bep[Proof of the {lower bound} of Theorem \ref{theorem_main_upper_tails}]
First, fix two configurations $(\vp, \psi) \in B$ such that $|A(\vp, \psi)| \le \inf_{(\vp',\psi')\in B} |A(\vp',\psi')| + \eps$. Because  $\inf_{(\vp',\psi')\in B} |A(\vp',\psi')|<\ff$ was assumed at the start of this section, also $|A(\varphi,\psi)|<\ff$ has to be satisfied. Now, let $\delta>0$ such that a) and b) from Lemma \ref{lemma_approximation_optimal_configurations} are satisfied. Under the event $\cap_{x\in\tau_n\vp} \{X(B_1(x))=1\}$, we can find $(z_x)_{x\in\tau_n\psi} \subseteq B_1(0)$ such that $\{x+z_x\} = X\cap B_1(x)$ for each $x\in\tau_n\psi$. Further, if $n$ is so large that $\t_n \de \ge 1$, under $\{X(\t_nA_{\delta}^-)=0\}$ it is guaranteed by \eqref{INF} that for each $x+z_x\in\{y+z_y \co y\in\tau_n\vp\}\cup\cup_{w\in\tau_n\vp} (\mc E_{w+z_w}(\{y+z_y \co y\in\tau_n\psi\}))$
\begin{equation}\label{inequality_lower_bound_subset}
\mc E(X-(x+z_x)) \supseteq \mc E(\{y+z_y \co y\in\tau_n\psi\}-(x+z_x)).
\end{equation}
Then, also if $\t_n \de \ge 1$, Lemma \ref{lemma_approximation_optimal_configurations} b) and \eqref{inequality_lower_bound_subset} give that
\begin{equation}\label{inequality_lower_bound_functional}
\sum_{x\in\tau_n\psi} \xi^{(\a)}(X-(x+z_x)) \ge \sum_{x\in\tau_n\vp} \xi^{(\a)}(\{y+z_y \co y\in\tau_n\psi\}-(x+z_x)) > \tau_n^\a/(1+\eps) = (r+\eps)n^d.
\end{equation}
Note that the index set in the sum before the first inequality in \eqref{inequality_lower_bound_functional} contains more points than the one after it. The reason for this is that when adding a point outside of the influence zone of $(\tau_n\varphi,\tau_n\psi)$, our framework for graphs does not exclude new edges from being created between two already existing nodes in $\tau_n\psi$. While it is admittedly hard to come up with an actual example for a graph for which the following is possible, it might potentially happen that when adding the points from $X\setminus(\tau_n\psi)$ to $\tau_n\psi$,  an additional edge arises from a point in $\tau_n\psi\setminus(\tau_n\vp)$ to a point in $\tau_n\vp$. In an undirected graph this could have the effect that the power-weighted edge lengths of some edges outgoing from points in $\tau_n\vp$ are only taken into account with the factor $1/2$ on the left-hand side of  the first inequality in \eqref{inequality_lower_bound_functional}, while being considered with their full weight on the right-hand side of it. Summing over all points in $\tau_n\psi$ avoids this issue.

As the remark after Lemma \ref{lemma_diam_influence_zone} suggets, we can assume that all of the occurred sets and configurations are contained in $U_n$ from which point \eqref{inequality_lower_bound_functional} implies that $G_{2,n}(\delta) \subseteq \{H_n(U_n,\R^d) > r + \eps\}.$

We now  define the set $W_n^- \coloneqq \big[(\log n)^2, b_n - (\log n)^2\big]$
and assert that under the event $G_{1, n}$ from \eqref{definition_G1}, we have 
\begin{align}
\label{wnm_eq}
H_n(W_n\sm W^-_n, W_n) < \e/2
\end{align}
for all large $n$. Once \eqref{wnm_eq} is established, we can conclude the proof of the lower bound of Theorem \ref{theorem_main_upper_tails}. Indeed, under $E_n^\ms{good}$ each box in $W_n\setminus W_n^-$ contains at least $c_\ms{max}$ Poisson points and therefore, if $n$ is chosen large, each of the cones around an $x\in X\cap W_n^-$ has to contain $c_\ms{max}$ Poisson points before the base of the cone leaves $W_n$, which more formally means that $\cup_{i \le I_d} \big ((S_i+x)\cap B_{\SS_i (X-x)}  (x)\big) \subseteq W_n$. Thus, under $E_n^\ms{good}$, again due to \eqref{STA}, we get that
$\xi^{(\a)}(X\cap W_n-x) = \xi^{(\a)}(X-x)$ for all points $x\in X\cap W^-_n$ if $n$ is sufficiently large. In other words, the layer of boxes containing points would not admit the score of points in $W_n^-$ being influenced by any points outside of $W_n$. With \eqref{wnm_eq} we get that under $G_{1,n}$
\begin{align*}
H_n(W_n,\R^d) \ge H_n(W^-_n,\R^d) = H_n(W^-_n, W_n) = H_n(W_n, W_n) - H_n(W_n\sm W^-_n, W_n) 
>  \mu_\a - \eps.
\end{align*}
Besides that, $G_{1, n}$ and $G_{2,n}(\delta)$ are independent for large $n$. Next, shifting the coordinate system shows that
$\P(H_n > \mu_\a + r) = \P(H_n([0,n]^d,\R^d) > \mu_\a + r).$
Hence,
\begin{align*}
\P(H_n > \mu_\a + r) &\ge \P\big(H_n(U_n,\R^d) > r + \eps, H_n(W_n,\R^d) > \mu_\a - \eps\big) \ge \P(G_{2,n}(\delta), G_{1,n}) = \P(G_{2,n}(\delta)) \P(G_{1,n}).
\end{align*}
Using Lemmas \ref{lemma_weak_law_of_large_numbers} and \ref{lemma_probability_points_responsible_for_excess_weight}, it follows that
\begin{equation}\label{inequality_lower_bound}
\liminf_{n\to\infty} \frac1{n^{d^2/\a}} \log \P(H_n > \mu_\a + r) \ge -\big(\inf_{(\vp,\psi)\in B} |A(\vp,\psi)|+\eps\big)(r+\eps)^{d/\a}
\end{equation}
and letting $\eps\downarrow 0$ gives the asserted result.

It remains to prove \eqref{wnm_eq} under the event $G_{1, n}$. To that end, we recall that 
$$H_n(W_n\sm W_n^-,W_n) = \frac1{n^d} \sum_{x\in X\cap (W_n\sm W^-_n)} \xi^{(\a)}(X\cap W_n - x).$$
Henceforth, we bound the summands on the right-hand side separately in the cases where $\dist(x,\partial W_n)\ge c\log n$ and where $\dist(x,\partial W_n)< c\log n$ for a suitable $c > 0$.

First, consider the case  $\dist(x,\partial W_n)\ge c\log n$. If we cut off the cone $S_i + x$ at a distance $c\log n$ for large enough $c$, then it still contains one of the boxes $Q_j'$. By definition of the event $E_n^\ms{good}$ each of these boxes contains at least $c_\ms{max}$ nodes. Therefore, $\mc S_i(X\cap W_n-x)$ is of order $\log n$ for all $1\le i\le I_d$.

Now, consider the case that $\dist(x,\partial W_n) < c\log n$.  If a cone that arises from $x$ does not intersect $W_n$ anymore after a distance from the apex of order $\log n$, then it contains Poisson points of $X\cap W_n$ only up until a distance of order $\log n$. Since by \eqref{STA}, none of the lateral boundaries of any cone are parallel to an axis of the coordinate system, we obtain that otherwise 
the cone envelopes a whole box $Q_j'$ after a distance from the apex of order $\log n$. Then, $\mc S_i(X\cap W_n-x)$ is of order $\log n$ as argued for above. Hence, after choosing $n$ sufficiently large, \eqref{STA} yields a finite configuration $\theta_x$ with $\theta_x\subseteq B_{(\log n)^2}(x)\cap W_n$ satisfying that $\mc E(X\cap W_n-x)=\mc E(\theta_x-x)$. Together with the bounded node degree, we have for $n$ large
\begin{align*}
	n^dH_n(W_n\sm W_n^-,W_n) \le \hspace{-.5cm}  \sum_{\substack{x\in X\cap (W_n\sm W^-_n)\\y\in\mc E(\theta_x -x)}} \hspace{-.5cm}|y| 
	&\le X(W_n\sm W^-_n) c_\ms{max} (\log n)^{2\a} 
	\le  (\log n)^{2d}\left\lceil\frac{n^{d-1}}{(\log n)^{d-3}}\right\rceil c_\ms{max} (\log n)^{2\a},
\end{align*}
where the final inequality follows from \eqref{number_of_boxes_to_cover_cube}, the upper bound on the number of Poisson points in each box $Q_i'$. Hence, we can choose $n$ sufficiently large to assure that $H_n(W_n\sm W_n^-,W_n)<\eps/2$.

\enp

The key ingredient to prove Lemma \ref{lemma_weak_law_of_large_numbers} is a weak law of large numbers for Poisson functionals from \cite{weakLLN}.

\bep[Proof of Lemma \ref{lemma_weak_law_of_large_numbers}]
We separately consider each of the two events that form $G_{1,n}$ when intersecting. For $E_n^{\ms{good}}$ we use a Poisson bound from \cite[Lemma 1.2]{poisson_conc} and calculate for $n$ sufficiently large
\begin{align}\label{convergence_regular_boundary_E1}
\begin{split}
\P(E_n^\ms{good}) &= 1-\P\big(X(Q'_i) < c_\ms{max} \text{ or } X(Q'_i) \ge (\log n)^{2d} \text{ for some } i\big) \\
	&\ge 1-\left\lceil\frac{n^{d-1}}{(\log n)^{d-3}}\right\rceil \bigg(e^{-(\log n)^d} \sum_{i \le c_\ms{max}-1} \frac{(2\log n)^{di}}{i!} + e^{-\frac12(\log n)^{2d}\log\big(\frac{(\log n)^{2d}}{(2\log n)^d}\big)}\bigg) \\
&\ge 1-n^{d-1}\big(e^{-(\log n)^d} (\log n)^{d c_\ms{max}+2} + e^{-\frac12(\log n)^{2d}}\big).
\end{split}
\end{align}
%The same calculation works for $E_n^\text{good,2}$ and, hence, for $n$ large
%\begin{equation}\label{convergence_regular_boundary_E2}
%\P(E_n^\text{good,2}) \ge 1-n^{d-1}\big(e^{-(\log n)^d} (\log n)^{d(k+1)+2} + e^{-\frac12(\log n)^{2d}}\big)
%\end{equation}
Next, we deal with $\{H_n(W_n,W_n) > \mu_\a - \eps/2\}$. Under the condition that $X(W_n) = b_n^d$,
we can deduce
\begin{align*}
\big(H_n(W_n,W_n)\ \big|\ X(W_n) = b_n^d\big) &\overset d= \frac{b_n^d}{n^d} \frac1{b_n^d} \sum_{i \le b_n^d} \xi^{(\a)}\big(\{X_1^{(n)},\dots,X_{b_n^d}^{(n)}\} - X_i^{(n)}\big) \\
&\overset d= \frac{b_n^d}{n^d} \frac1{b_n^d} \sum_{i \le b_n^d} \xi^{(\a)}\big(b_n\{\wh X_1,\dots,\wh X_{b_n^d}\} - b_n\wh X_i\big)
\end{align*}
for $X_1^{(n)}, X_2^{(n)},\dots$ i.i.d.\ uniform random variables on $W_n$ and $\wh X_1,\wh X_2,\dots$ i.i.d.\ uniform random variables on $[0,1]^d$. In order to apply \cite[Theorem 2.1]{weakLLN}, we need to check the moment condition, i.e., that for $p>2$
\begin{equation}\label{inequality_moment_condition}
\sup_{n\ge 1} \E\big[\xi^{(\a)}(b_n\{\wh X_1,\dots, \wh X_{b_n^d}\}- b_n\wh X_1)^p\big] <\infty.
\end{equation}
We can use the bound on the node degree to get
\begin{align}\label{inequality_out_neighbor}
\begin{split}
&\E\big[\xi^{(\a)}(b_n\{\wh X_1,\dots, \wh X_{b_n^d}\}- b_n\wh X_1)^p\big] = \int_0^\infty \P\big(\xi^{(\a)}(b_n\{\wh X_1,\dots, \wh X_{b_n^d}\}- b_n\wh X_1)^p > s\big)ds \\
&\le \int_0^\infty \P\big(|y| > (s^{1/p}/c_\ms{max})^{1/\a} \text{ for some }y\in \mc E(b_n\{\wh X_1,\dots, \wh X_{b_n^d}\}- b_n\wh X_1)\big)ds.
\end{split}
\end{align}
Next, from \eqref{STA} we can deduce that, for every $s>0$, if $b_n \wh X_1$ has an out-neighbor among $b_n\{\wh X_2,\dots, \wh X_{b_n^d}\}$ that is farther away than $(s^{1/p}/c_\ms{max})^{1/\a}$, one of the cones arising from $b_n \wh X_1$ has to extend until at least a distance of $(s^{1/p}/c_\ms{max})^{1/\a}$ from its apex before it contains $c_\ms{STA}$ vertices. More precisely, there has to be an $i\le I_d$ such that $\mc S_i(b_n\{\wh X_1,\dots, \wh X_{b_n^d}\}- b_n\wh X_1) > (s^{1/p}/c_\ms{max})^{1/\a}$. Additionally, the intersection of $W_n$ and $(S_i + b_n \wh X_1)\setminus B_{(s^{1/p}/c_\ms{max})^{1/\a}}(b_n \wh X_1)$ cannot be empty since the mentioned out-neighbor has to be within $W_n$. Therefore, under the event from the last line of \eqref{inequality_out_neighbor}, it is implied by the definition of $\mc S_i(\cdot)$ that for some $i \le I_d$ it holds that $b_n\{\wh X_1,\dots, \wh X_{b_n^d}\}\cap(S_i+b_n\wh X_1)\cap B_{(s^{1/p}/c_\ms{max})^{1/\a}}(b_n\wh X_1)$ contains at most $c_\ms{max}$ points, while $W_n\cap (S_i+b_n\wh X_1)\sm B_{(s^{1/p}/c_\ms{max})^{1/\a}}(b_n\wh X_1)\neq\emptyset$. With these arguments we arrive at
\begin{align}\label{inequality_far_neighbor_to_cones}
\begin{split}
&\int_0^\infty \P\big(|y| > (s^{1/p}/c_\ms{max})^{1/\a} \text{ for some }y\in \mc E(b_n\{\wh X_1,\dots, \wh X_{b_n^d}\}- b_n\wh X_1)\big)ds \\
&\le \sum_{i \le I_d} \int_0^\infty \P\big(\#(b_n\{\wh X_1,\dots, \wh X_{b_n^d}\}\cap(S_i+b_n\wh X_1)\cap B_{(s^{1/p}/c_\ms{max})^{1/\a}}(b_n\wh X_1)) \le c_\ms{max} \\
&\qquad\qquad\qquad\ \text{and } W_n\cap (S_i+b_n\wh X_1)\sm B_{(s^{1/p}/c_\ms{max})^{1/\a}}(b_n\wh X_1) \ne\emptyset \big)ds.
\end{split}
\end{align}
Further, since the cones do not have lateral boundaries parallel to any axes, under the event after the last inequality in \eqref{inequality_far_neighbor_to_cones}, the volume of the set $W_n\cap (S_i+b_n\wh X_1)\cap B_{(s^{1/p}/c_\ms{max})^{1/\a}}(b_n\wh X_1)$ is of order $s^{d/(p\a)}$ and therefore at least $c s^{d/(p\a)}$ for all $i$, where $c>0$ depends only on the layout of the cones, $\alpha$ and $c_\ms{max}$. Thus, due to the independence of $\wh X_2,\dots,\wh X_{b_n^d}$, when conditioned on $\wh X_1$, we can bound the probability of the event after the last inequality in \eqref{inequality_far_neighbor_to_cones} by the probability of a binomial random variable consisting of $b_n^d-1$ trials with success probability
$$\P(b_n \wh X_2 \in (S_i+b_n\wh X_1)\cap B_{(s^{1/p}/c_\ms{max})^{1/\a}}(b_n\wh X_1) \mid\, W_n\cap (S_i+b_n\wh X_1)\sm B_{(s^{1/p}/c_\ms{max})^{1/\a}}(b_n\wh X_1) \ne\emptyset) \ge cs^{d/(p\a)}/b_n^d$$
realizing a value of at most $c_\ms{max}-1$. Thus, by a binomial concentration bound \cite[Lemma 1.1]{poisson_conc}
\begin{align*}
&\E\big[\xi^{(\a)}(b_n\{\wh X_1,\dots, \wh X_{b_n^d}\}- b_n\wh X_1)^p\big] \\
&\le I_d \int_0^\infty \P\big(\ms{Bin}(b_n^d-1, c s^{d/(p\a)} / b_n^d) \le c_\ms{max}-1\big)ds \le I_d \int_0^\infty \P\big(\ms{Bin}(b_n^d, c s^{d/(p\a)} / b_n^d) \le c_\ms{max}\big)ds \\
&\le I_d (c_\ms{max}/c)^{p\a/d} + I_d \int_{(c_\ms{max}/c)^{p\a/d}}^\infty \P\big(\text{Bin}(b_n^d, c s^{d/(p\a)} / b_n^d) \le c_\ms{max}\big)ds \\
&\le I_d (c_\ms{max}/c)^{p\a/d} +I_d \int_{(c_\ms{max}/c)^{p\a/d}}^\infty \exp\big(-c s^{d/(p\a)} (1- \tfrac{c_\ms{max}}{c s^{d/(p\a)}} + \tfrac{c_\ms{max}}{c s^{d/(p\a)}} \log(\tfrac{c_\ms{max}}{c s^{d/(p\a)}})\big)ds < \infty.
\end{align*}
In particular, the bound does not depend on $n$. Therefore, the moment condition \eqref{inequality_moment_condition} is satisfied. Now, \cite[Theorem 2.1]{weakLLN} gives
$$\frac1{b_n^d} \sum_{i=1}^{b_n^d} \xi^{(\a)}\big(b_n\{\wh X_1,\dots,\wh X_{b_n^d}\} - b_n\wh X_i\big) \overset{P}{\longrightarrow} \mu_\a$$
and since $b_n^d/n^d \overset{n\to\infty}{\longrightarrow} 1$, it follows that
\begin{equation}\label{convergence_weak_law}
\lim_{n\to\infty} \P\big(H_n(W_n,W_n) > \mu_\a - \eps/2\, |\, X(W_n) = b_n^d\big) = 1.
\end{equation}
Now, we can use the union bound to arrive at
\begin{align*}
\P(G_{1,n}) &\ge \P\big(H_n(W_n,W_n) > \mu_\a - \eps/2\big) + \P(E_n^\ms{good})  - 1 \\
&\ge \P\big(H_n(W_n,W_n) > \mu_\a - \eps/2\, |\, X(W_n) = b_n^d\big) \P(X(W_n) = b_n^d) + \P(E_n^\ms{good})  - 1\\
&\ge \P\big(H_n(W_n,W_n) > \mu_\a - \eps/2\, |\, X(W_n) = b_n^d\big) \frac{\exp\left(-\frac1{12 b_n^d}\right)}{\sqrt{2\pi b_n^d}} + \P(E_n^\ms{good})  - 1,
\end{align*}
where in the last line, we used \cite[Lemma 1.3]{poisson_conc}.
Due to \eqref{convergence_weak_law} we can assume $n$ large enough so that
$$\P\big((H_n(W_n,W_n) > \mu_\a - \eps/2\, |\, X(W_n) = b_n^d\big) \ge 1/2.$$
Hence, together with  \eqref{convergence_regular_boundary_E1}, we get that
$$\P(G_{1,n}) \ge \frac12\frac{\exp\left(-\frac1{12b_n^d}\right)}{\sqrt{2\pi b_n^d}} - n^{d-1}\big(e^{-(\log n)^d} (\log n)^{d c_\ms{max}+2} + e^{-\frac12(\log n)^{2d}}\big),$$
as asserted.
\enp

We continue with the proof of Lemma \ref{lemma_approximation_optimal_configurations}.

\bep[Proof of Lemma \ref{lemma_approximation_optimal_configurations}]
First, we show that
\begin{equation}\label{equality_approximation_configuration}
\bigcap_{\delta>0}\underbrace{\bigcup_{(z_y)_{y\in\psi}\subseteq B_\delta(0)} A(\{x+z_x \colon x\in\vp\},\{x+z_x \colon x\in\psi\})}_{\eqqcolon A_\delta(\vp,\psi)} \subseteq A(\vp,\psi) \cup D(\psi)
\end{equation}
where we recall that $D(\psi) = \{y \in \R^d \co \vp \cup\{y\} \in N_{\#\vp +1}\}$. The subset relation in \eqref{equality_approximation_configuration} holds because if we let $y\in\cap_{\delta>0} A_\delta(\vp,\psi)\sm D(\psi)$, then for every $\delta\in(0,1)$ there exists a family $(z_w)_{w\in\psi}\subseteq B_\delta(0)$ such that $y\in A(\{w+z_w \co w\in\vp\},\{w+z_w \co w\in\psi\})$. Hence,
$$\mc E_{x+z_x}(\{w+z_w \co w\in\psi\}) \not\subseteq \mc E_{x+z_x}(\{w+z_w \co w\in\psi\}\cup\{y\})$$
for some $x+z_x\in\{w+z_w \co w\in\vp\}\cup\cup_{v+z_v\in\{w+z_w \co w\in\vp\}} (\mc E_{v+z_v}(\{w+z_w \co w\in\psi\}))$. Since $y\not\in D(\psi)$, we can apply \eqref{CON} to both sides, choosing $\delta$ sufficiently small, which gives $\mc E_x(\psi) \not\subseteq \mc E_x(\psi\cup\{y\})$ for some $x\in\vp\cup\cup_{v\in\vp} (\mc E_v(\psi))$ and therefore $y\in A(\vp,\psi)$. Since $|D(\psi)|=0$, we deduce from \eqref{equality_approximation_configuration} that for $\delta$ sufficiently small we have
$
|A_\delta(\vp,\psi)|\leq |A(\vp,\psi)| + \eps.
$

To prove part b) note that due to $\psi$ being finite, we can use \eqref{CON} and find $\delta\in(0,1)$ small enough such that for all choices of $(z_x)_{x\in\psi} \subseteq B_\delta(0)$ we have
\begin{equation}\label{equality_continuous_application}
\{w - z_w \co w\in\mc E_{x+z_x}(\{y+z_y \co y\in\psi\})\} = \mc E_x(\psi).
\end{equation}
This means the graph looks the same despite some small noise of at most $\delta$ for every node. But the finiteness of the configuration combined with \eqref{equality_continuous_application} guarantees that for $\delta$ sufficiently small
$$\inf_{(z_x)_{x\in\psi} \subseteq B_\delta(0)} \sum_{x\in\vp} \xi^{(\a)}(\{y+z_y \colon y\in\psi\}-(x+z_x)) > 1/(1+\eps).$$
\enp

What follows is the proof of Lemma \ref{lemma_diam_influence_zone}.
\bep[Proof of Lemma \ref{lemma_diam_influence_zone}]
Let $(\varphi,\psi)\in B$ be such that $|A(\varphi,\psi)|<\infty$. The key step is to  construct a finite set of points $\th \su \R^d$ and a scalar $R > 0$ such that for all $x\in\eta\coloneqq\varphi \cup\cup_{z\in\varphi} (\EE_z (\psi))$ and $(z_y)_{y\in\psi}\subseteq B_\de(0)$ we have i) $\RR\big((\{y+z_y\co y\in\psi\} \cup \th) - (x + z_x)\big) \le R$, and ii) 
\begin{align}
	\label{th_eq}
	\EE\big(\{y+z_y\co y\in\psi\}-(x+z_x)\big)\subseteq\EE\big((\{y+z_y\co y\in\psi\}\cup\th )-(x+z_x)\big).
\end{align}
Once $\th$ is constructed, we assert that $A_\de(\vp, \psi) \su \bigcup_{x \in \eta} B_{R + 1}(x)$. Indeed, for any $v \in \R^d\sm \bigcup_{x \in \eta} B_{R + 1}(x)$, the definition of the stabilization radius implies that 
\begin{equation}\label{subsetequality_diameter}
	\EE\big((\{y+z_y\co y\in\psi\}\cup\theta)-(x+z_x)\big) = \EE\big((\{y+z_y\co y\in\psi\}\cup\theta\cup\{v\})-(x+z_x)\big).
\end{equation}
Hence, combining \eqref{th_eq}, \eqref{subsetequality_diameter} and \eqref{STA} gives that
$$	\EE\big(\{y+z_y\co y\in\psi\}-(x+z_x)\big) \su \EE\big((\{y+z_y\co y\in\psi\}\cup\{v\})-(x+z_x)\big),$$
thereby showing the asserted $v \not\in A_\de(\vp, \psi)$.

It remains to prove the existence of $R > 0$ and $\th \su \R^d$. To that end, first note that $|T_i(x,\de)| = \infty$ for all $i \le I_d$ and $x\in\R^d$, where $T_i(x,\de)\coloneqq \cap_{y\in B_\de(x)} (S_i+y)$. Since $|A_\de(\varphi,\psi)|<\infty$, due to Lemma \ref{lemma_approximation_optimal_configurations} a) if $\de$ is chosen appropriately, it follows that for every $x\in \eta$ and $i\le I_d$, there are distinct $w_{x,i}^{(1)},\dots,w_{x,i}^{(c_\ms{STA})}\in T_i(x,\de)\setminus A_\de(\varphi,\psi)$. Then, defining $\theta\coloneqq \{w_{x,i}^{(j)} \co x\in\eta, i \le I_d, j \le c_\ms{STA}\}$, we note that \eqref{INF} implies property \eqref{th_eq}. Now, set
$$R\coloneqq \max_{\substack{x\in\eta, \, z\in B_\de(x)\\  i\le I_d, j\le c_\ms{STA}}} c_\ms{STA}|z-w_{x,i}^{(j)}| \le \max_{\substack{x\in\eta\\  i \le I_d,  j\le c_\ms{STA}}} c_\ms{STA} |x-w_{x,i}^{(j)}|+c_\ms{STA} \de,$$
and note that the finiteness of the configurations in $B$ implies that $R < \ff$. Then, the definition of the stabilization radius yields $\RR\big((\{y+z_y\co y\in\psi\} \cup \th) - (x + z_x)\big) \le R$, 
as asserted.

\enp

Finally, we show Lemma \ref{lemma_probability_points_responsible_for_excess_weight}.

\bep[Proof of Lemma \ref{lemma_probability_points_responsible_for_excess_weight}]
Let $(\vp,\psi)\in B$ and $\delta\in(0,1)$ be given according to the setting of Lemma \ref{lemma_approximation_optimal_configurations}. First note that $B_1(x) \su B_{\tau_n \de}(x)$ for sufficiently large $n\ge 1$. Moreover, the events $\{X(B_1(x)) = 1 \text{ for all } x\in\tau_n\psi\}$ and $\{X(A_{\delta, n}^-) = 0\}$ are independent.  Thus, we can examine them separately and start with the first one. We assume that $n$ is chosen large enough such that all of the balls around points in $\tau_n\psi$ are disjoint. Therefore,
$$\P\big(X(B_1 (x)) = 1 \text{ for all } x\in\tau_n\psi\big) = \prod_{x\in\tau_n\psi}\P\big(X(B_1 (x)) = 1\big) = \kappa_d  ^{\#\psi} e^{-\#\psi \kappa_d}.$$
For the second event, Lemma \ref{lemma_approximation_optimal_configurations} a) yields
$$\P\big(X(A_{\delta, n}^-) = 0\big) \ge \P\big(X(\t_nA_{\delta}(\vp,\psi)) = 0\big) \ge \exp\big(-\big( |A(\vp,\psi)|+\eps\big)\tau_n^d\big).$$
All of this combined shows that for large enough $n$
\begin{align*}
	\frac1{n^{d^2/\a}} \log \P(G_{2,n}(\delta)) &\ge \frac1{n^{d^2/\a}} \log\big( \kappa_d^{\#\psi} e^{-\#\psi \kappa_d}\big) - \big(\inf_{(\vp,\psi)\in B} |A(\vp,\psi)|+\eps\big)((r+\eps)(1+\eps))^{d/\a} \\
&\overset{n\to\infty}{\longrightarrow} - \big(\inf_{(\vp,\psi)\in B} |A(\vp,\psi)|+\eps\big)((r+\eps)(1+\eps))^{d/\a},
\end{align*}
as asserted.
\enp

%
%PART B
%
\section{Proof of Theorem \ref{theorem_conditioned_convergence}}\label{sec_proof_2}
The proof of  Theorem \ref{theorem_conditioned_convergence} consists mainly of a refinement of the steps from the proofs of Theorem \ref{theorem_main_upper_tails} and Lemma \ref{lemma_sum_optimization_problem}.

\bep[Proof of Theorem \ref{theorem_conditioned_convergence} a)]
To begin with, let $\delta\in (0,1)$ and for $\eps\in(0,(1-d/\a)/(2\a))$, set $h_{n,\eps} \coloneqq \lfloor n^{d^2/\a-\eps}\rfloor$.  We  put $H_n' \coloneqq n^{-d}\sum_{i \le h_{n,\eps}} Z_n^{(i)}$ and
will separately look at the numerator and denominator of
$$
\P\big(H_n' < r(1-\delta)\, \big|\, H_n > \mu_\alpha + r\big) = \frac{\P(H_n' < r(1-\delta), H_n > \mu_\alpha + r)}{\P(H_n > \mu_\alpha + r)}
$$
and prove that this ratio tends to zero. To start with the numerator, recall how Lemma \ref{lemma_sum_optimization_problem} was used in the proof of Theorem \ref{theorem_main_upper_tails}. The event $\{H_n>\mu_\a+r\}$ was split up in small and large contributions. Instead of $\eps$, here we use an arbitrary $\tilde\eps < (1-d/\a)/(2\a) \wedge \delta$ to divide the term $\mu_\a n^d+rn^d$. Then, since $\tilde\eps<\delta$, we have
$$\big\{H_n' < r(1-\delta)\big\} \cap \bigg\{\sum_{x\in\mathcal{J}_n^{(a)}(X)} \xi^{(\a)}(X-x) \ge rn^d(1-\tilde\eps), J_n^{(a)}(X) \le h_{n,\eps}\bigg\} = \emptyset.$$
 Using this, Lemma \ref{lemma_hoeffding} and Lemma \ref{lemma_poisson_conc} similarly as in the proof of Theorem \ref{theorem_main_upper_tails} gives that
\begin{equation}\label{inequality_sum_with_max_bound}
\limsup_{n\to\infty} \frac1{n^{d^2/\a}} \log\P\big(H_n' < r(1-\delta),H_n > \mu_\alpha + r\big) =-\infty.
\end{equation}
For the denominator, after additionally assuming that $\tilde\eps < \mu_\a$, we deduce from \eqref{inequality_lower_bound} that
\begin{equation}\label{inequality_sum_below}
	\P(H_n > \mu_\a + r) \geq \exp(-\g(\tilde \e)n^{d^2/\a} + o(n^{d^2/\a})),
\end{equation}
where $\g(\tilde \e) \coloneqq (\inf_{(\varphi,\psi)\in B} |A(\varphi,\psi)|+\tilde\eps\big)(r+\tilde\eps)^{d/\a}$. Together, \eqref{inequality_sum_with_max_bound} and \eqref{inequality_sum_below} imply that for any $c>0$, if $n$ is chosen large enough, we have that
\begin{align*}
	\P\big(H_n' < r(1-\delta)\, \big|\, H_n > \mu_\alpha + r\big) &\leq \frac{\exp(-cn^{d^2/\a} + o(n^{d^2/\a}))}{\exp(-\g(\tilde \e)n^{d^2/\a} + o(n^{d^2/\a}))},
\end{align*}
which indeed converges to $0$ when letting $n$ go to infinity, provided that $c > \g(\tilde\e)$.

Next, we prove that the statement about the other side holds, i.e., 
$$\P\big(H_n' > r(1+\delta)\, \big|\, H_n > \mu_\alpha + r\big) \overset{n\to\infty}{\longrightarrow} 0$$
To that end, we note that the proof of Lemma \ref{lemma_sum_optimization_problem} extends without any changes to the case where we replace $\mc J_n^{(a)}(X)$ by the set of the nodes with the $h_{n, \e}$ largest scores. Then, applying this result with  $\t =rn^d(1+\delta) $ and $m =h_{n,\eps}$ yields
\begin{align*}
&\P\big(H_n' > r(1+\delta), \RR_{3n} \le n\big) \\
&\leq ((I_d c_\ms{max} + 1)^2 h_{n,\eps})^2 (5n)^{d 2 (I_d c_\ms{max} + 1)^2 h_{n,\eps}}\exp\Big(- (rn^d(1+\delta))^{d/\a} \inf_{(\varphi,\psi)\in B} |A(\varphi,\psi)|\Big) \\
&= \exp\Big(- (rn^d(1+\delta))^{d/\a} \inf_{(\varphi,\psi)\in B} |A(\varphi,\psi)| + o(n^{d^2/\a})\Big).
\end{align*}
Proceeding similar to the proof of Theorem \ref{theorem_main_upper_tails}, we get the bound for the numerator and can estimate
$$\P\big(H_n' > r(1+\delta)\, \big|\, H_n > \mu_\alpha + r\big) \leq \frac{\exp\big(- (rn^d(1+\delta))^{d/\a} \inf_{(\varphi,\psi)\in B} |A(\varphi,\psi)| + o(n^{d^2/\a})\big)}{\exp(-\g(\tilde\eps)n^{d^2/\a} + o(n^{d^2/\a}))}.$$
Now, choosing $\tilde\eps$ small enough to achieve that the bound for the numerator converges to $0$ faster than the bound for the denominator gives  the claimed convergence. Thus,
$$\P\big(|{H_n'}/r - 1| > \delta\, \big|\, H_n > \mu_\alpha + r\big) \overset{n\to\ff}{\longrightarrow} 0,$$
as asserted.
\enp

Showing part b) mainly requires redoing the steps of part a). Nevertheless, it is a bit more challenging since we need to replicate Lemma \ref{lemma_sum_optimization_problem} in a slightly extended form that incorporates the additional bound for the sum of the largest scores within the sample space.

\bep[Proof of Theorem \ref{theorem_conditioned_convergence} b)]
Let $m_0>0$ satisfy \eqref{inequality_condensation_uniqueness}. Let $\delta, \eps > 0$ be chosen as in the proof of part a). This time, let $\tilde\eps < (1-d/\a)/(2\a)\wedge \delta/2$ be arbitrary and define $H_n' \coloneqq n^{-d}\sum_{i=1}^{m_0}  Z_n^{(i)}$ and $H'(\vp, \psi) \coloneqq \sum_{i=1}^{m_0}  Z^{(i)}(\vp, \psi)$.  As in the proof of Lemma \ref{lemma_sum_optimization_problem}, we can show that
\begin{align}\label{inequality_sum_optimization_incl_max_bound}
\begin{split}
&\P\bigg(H_n' < (1-\tilde\eps)r\tfrac{1-\delta}{1-\delta/2}, \sum_{x\in\mathcal{J}_n^{(a)}(X)} \xi^{(\a)}(X-x) \geq (1-\tilde\eps)rn^d, J_n^{(a)}\le \lfloor n^{d^2/\a-\eps}\rfloor, \RR_{3n}(X) \leq n\bigg) \\
&\leq k_{n,\eps} \exp\Big(-((1-\tilde\eps)rn^d)^{d/\a}\inf_{(\varphi,\psi)\in B,\, H'(\varphi,\psi) < (1-\delta)/(1-\delta/2)} |A(\varphi,\psi)|\Big),
\end{split}
\end{align}
where $k_{n,\eps} \coloneqq (I_d c_\ms{max} + 1)^4 \lfloor n^{d^2/\a-\eps}\rfloor^2 (5n)^{d 2 (I_d c_\ms{max} + 1)^2 \lfloor n^{d^2/\a-\eps}\rfloor} \in e^{o(n^{d^2/\a})}$. Further, repeating the arguments from the proof of Theorem \ref{theorem_main_upper_tails} as we did to get (\ref{inequality_sum_with_max_bound}), but replacing Lemma \ref{lemma_sum_optimization_problem} with \eqref{inequality_sum_optimization_incl_max_bound}, we arrive at
\begin{align}\label{inequality_sum_with_max_bound2}
\begin{split}
\P\big(H_n' < r(1-\delta), H_n > \mu_\alpha + r\big) &\le \P\big(H_n' < (1-\tilde\eps)r\tfrac{1-\delta}{1-\delta/2}, H_n > \mu_\alpha + r \big) \\
&\le \exp\Big(-((1-\tilde\eps) rn^d)^{d/\a} \hspace*{-0.4725cm} \inf_{(\varphi,\psi)\in B,\, H'(\varphi, \psi) < (1-\delta)/(1-\delta/2)} |A(\varphi,\psi)| + o(n^{d^2/\a})\Big),
\end{split}
\end{align}
which is sufficient for dealing with the numerator.

For the denominator, we can reuse the inequality stated in \eqref{inequality_sum_below} with the assumption $\tilde\eps <\mu_\a$. Next, because of \eqref{inequality_condensation_uniqueness} applied to $\delta' = 1-(1-\delta)/(1-\delta/2)$, we can demand $\tilde\eps$ to be small enough to assure that
\begin{equation}\label{inequality_to_guarentee_convergence}
((1-\tilde\eps) r)^{d/\a} \inf_{(\varphi,\psi)\in B,\, H'(\varphi, \psi) < 1-\delta'} |A(\varphi,\psi)| > \Big(\inf_{(\varphi,\psi)\in B} |A(\varphi,\psi)|+\tilde\eps\Big)((r+\tilde\eps))^{d/\a}.
\end{equation}
We proceed by plugging (\ref{inequality_sum_below}) and (\ref{inequality_sum_with_max_bound2}) into the fraction that arises from the conditional probability and get
$$\P\big(H_n' < r(1-\delta)\, \big|\, H_n > \mu_\alpha + r\big) \le \frac{\exp\big(-(1-\tilde\eps)^{d/\a} r^{d/\a}n^{d^2/\a} \inf_{(\varphi,\psi)\in B,\, H'(\varphi, \psi) < 1-\delta'} |A(\varphi,\psi)|+o(n^{d^2/\a})\big)}{\exp\big(-(\inf_{(\varphi,\psi)\in B} |A(\varphi,\psi)|+\tilde\eps)(r+\tilde\eps)^{d/\a}n^{d^2/\a} + o(n^{d^2/\a})\big)},$$
which converges to $0$ due to the assumed relation of the coefficients in (\ref{inequality_to_guarentee_convergence}).

The assertion on the upper tails, i.e., $\P\big(H_n' > r(1+\delta)\, \big|\, H_n > \mu_\alpha + r\big) \overset{n\to\infty}{\longrightarrow} 0$, follows analogously to part a).
\enp

\section*{Acknowledgment}
We would like to thank the two anonymous referees for providing us with valuable comments and suggestions for the manuscript. Further, the authors would like to acknowledge the financial support of the CogniGron research center and the Ubbo Emmius Funds (Univ.~of Groningen).

\bibliographystyle{abbrv}
%\footnotesize
\bibliography{refs}

\end{document}